\documentclass[11pt]{article}
\usepackage{amsmath,amsfonts,amsthm,amscd,amssymb,graphicx}
\usepackage{url}
\usepackage[utf8]{inputenc}
\usepackage{amsbsy}
\usepackage{graphicx, color}
\usepackage[draft]{todonotes}

\usepackage{comment}
\usepackage[text={7.0in,9.5in},centering,includefoot,foot=0.6in]{geometry}
\usepackage{verbatim}
\usepackage{float}
\usepackage{subfloat}
\usepackage{subcaption}
\newtheorem{Lemma}{Lemma}[section]

\newtheorem{Remark}[Lemma]{Remark}
\newtheorem{Theorem}[Lemma]{Theorem}

 {\begin{trivlist}\item[]\textbf{Proof#1 }}%
 {\hspace*{\fill}$\rule{0.3\baselineskip}{0.35\baselineskip}$\end{trivlist}}

 {\begin{trivlist}\item[]\textbf{Acknowledgments }}{\end{trivlist}}

\makeatletter\@addtoreset{equation}{section}\makeatother

\def\Im{\mathop\mathrm{Im}\nolimits}    

\newcommand{\rmd}{\mathrm{d}}           
\newcommand{\rmi}{\mathrm{i}}           
\newcommand{\uu}{\mathbf{u}}

\title{Selection of quasi-stationary states in the stochastically forced Navier-Stokes equation on the torus}
\author{Margaret Beck\footnote{Department of Mathematics and Statistics, Boston University, Boston, MA 02215. E-mail: mabeck@bu.edu, cooper@bu.edu, kspiliop@math.bu.edu. M.B. was partially supported by National Science Foundation (NSF) grant DMS 1411460 and K.S. was partially supported by NSF DMS 1550918.}, Eric Cooper\footnotemark[1], Gabriel Lord\footnote{Department of Mathematics, Heriot-Watt University, Edinburgh, EH14 4AS, UK. Email: g.j.lord@hw.ac.uk}, Konstantinos Spiliopoulos\footnotemark[1]}
\date{\vspace{-5ex}}

\begin{document}

\maketitle

\begin{abstract}
The stochastically forced vorticity equation associated with the two dimensional incompressible Navier-Stokes equation on $D_\delta:=[0,2\pi\delta]\times [0,2\pi]$ is considered for $\delta\approx 1$, periodic boundary conditions, and viscocity $0<\nu\ll 1$. An explicit family of quasi-stationary states of the deterministic vorticity equation is known to play an important role in the long-time evolution of solutions both in the presence of and without noise. Recent results show the parameter $\delta$ plays a central role in selecting which of the quasi-stationary states is most important. In this paper, we aim to develop a finite dimensional model that captures this selection mechanism for the stochastic vorticity equation. This is done by projecting the vorticity equation in Fourier space onto a center manifold corresponding to the lowest eight Fourier modes. Through Monte Carlo simulation, the vorticity equation and the model are shown to be in agreement regarding key aspects of the long-time dynamics. Following this comparison, perturbation analysis is performed on the model via averaging and homogenization techniques to determine the leading order dynamics for statistics of interest for $\delta\approx1$.
\end{abstract}

\section{Introduction}\label{S:Intro}

Consider the 2D incompressible Navier-Stokes equation,
\begin{equation}
\begin{aligned}
\frac{\partial\uu}{\partial t}  = \nu \Delta \uu &- (\uu \cdot \nabla) \uu - \nabla p \label{E:2dns} \\
\nabla \cdot \uu &= 0,
\end{aligned}
\end{equation}
on the possibly asymmetric torus $(x,y) \in D_\delta := [0, 2\pi\delta] \times [0, 2\pi]$ with $\delta \approx 1$, periodic boundary conditions, and viscosity $0 < \nu \ll 1$. To obtain the equivalent vorticity formulation of the equation, take the curl of the vector field $\uu$ and set $\omega = (0,0,1) \cdot (\nabla \times \uu)$ to find
\begin{equation}\label{E:2dvort}
\partial_t \omega = \nu \Delta \omega - \uu\cdot\nabla\omega, \qquad \uu = \begin{pmatrix} \partial_y(- \Delta^{-1}) \\ -\partial_x(- \Delta^{-1})\end{pmatrix} \omega.
\end{equation}

The relation between $\uu$ and $\omega$ is known as the Biot-Savart law. The periodic boundary conditions insure $\int_{D_\delta} \omega = 0$, and therefore $\Delta^{-1} \omega$ is well-defined. 

Adding random forcing to the system allows to account for stochasticity/genericity in the system, see for example \cite{Novikov1965,BensousanTemam73,Glatt-Holtz}. In particular, we add a stochastic forcing term to \eqref{E:2dvort} to obtain the stochastic 2D vorticity equation,
\begin{equation}\label{E:2dvortS}
\partial_t \omega = \nu \Delta \omega - \uu\cdot\nabla\omega+\frac{\partial \mathcal{W}}{\partial t}, \qquad \uu = \begin{pmatrix} \partial_y(- \Delta^{-1}) \\ -\partial_x(- \Delta^{-1})\end{pmatrix} \omega.
\end{equation}
The noise is white in time, colored in space, and takes the form, for $\vec{k}=(k_1,k_2)\neq(0,0)$,
\begin{align}\label{E:noise}
\mathcal{W}(t,x,y)=\sqrt{2\nu}\sum_{\vec{k}\in \mathcal{K}\subset \mathbb{Z}^{2}\backslash\{(0,0)\}}\sigma_{\vec{k}}e^{\rmi (k_{1}x/\delta+k_{2}y)}\beta_{\vec{k}}(t),
\end{align}
with spatial correlation $\sigma_{\vec{k}}$ and $\mathcal{K}$ to be commented on below. Here $\beta(t)=\{\beta_{\vec{k}}(t)\}$ is a collection of i.i.d. Wiener processes.

Notice that with the noise (\ref{E:noise}), the equation (\ref{E:2dvortS}) is now stochastic.  To insure the random vorticity remains real valued for all times $t\geq0$, the following complex conjugacy conditions are imposed, $\bar{\sigma}_{\vec{k}}=\sigma_{-\vec{k}}$ and $\bar{\beta}_{\vec{k}}=\beta_{-\vec{k}}$. Additional assumptions are often placed on the noise coefficients, $\sigma_{\vec{k}}$, to insure certain smoothness properties of solutions. In particular, we assume that there exist fixed positive constants $C_0$ and $\alpha_0$ such that $|\sigma_{\vec{k}}|\leq C_0 e^{-\alpha_0|\vec{k}|^2}$ so that solutions will then be analytic in space \cite{Mattingly02}.  Since the boundary conditions force solutions of the deterministic equation to satisfy $\int_{D_\delta} \omega =0$, we choose $\sigma_{(0,0)}=0$ so this property is preserved. Note that if $\sigma_{\vec{k}}=0$ for all $\vec{k}\in\mathbb{Z}^2$ then \eqref{E:2dvortS} reduces to the deterministic vorticity equation.

Although an $L^2$ energy estimate shows solutions of \eqref{E:2dvort} have a time-asymptotic rest state of zero, certain quasi-stationary states, known as bars and dipoles, rapidly attract nearby solutions and correspond to transient structures that play a key role in the long-time evolution of solutions \cite{BeckCooperSpiliopoulos, BeckWayne13, BouchetSimonnet09, grenier, IbrahimMaekawaMasmoudi17, LinXu,weiz, Yin}. These quasi-stationary states are members of an explicit family of functions given by,
\begin{equation}\label{E:family}
\omega(x,y,t) = e^{-\frac{\nu}{\delta^2} t}[ a_1 \cos(x/\delta) + a_2 \sin (x/\delta)] + e^{-\nu  t}[a_3 \cos(y) + a_4 \sin (y) ].
\end{equation}

If $\delta=1$, then any member of this family is an exact solution to the deterministic vorticity equation. If $\delta \neq 1$, then (\ref{E:family}) remains a solution if and only if $a_1 = a_2 = 0$ or if $a_3 = a_4 = 0$. These members, which only depend on one spatial variable, are called bar states, and they are also known as unidirectional or Kolmogorov flow. The $x$- and $y$-bar states are members of this family given by
\[
\omega_{xbar}(x,t) = e^{-\frac{\nu}{\delta^2}t} \sin (x/\delta), \qquad \omega_{ybar}(y ,t) = e^{-\nu t} \sin y,
\]
or similarly with sine replaced by cosine. The associated velocity fields are given by
\[
\uu_{xbar}(x,t) = -\delta e^{-\frac{\nu}{\delta^2}t} \left(
\begin{array}{c}
0\\
\cos(x/\delta)\\
\end{array}
\right), \qquad \uu_{ybar}(y ,t) = e^{-\nu t} \left(
\begin{array}{c}
\cos y\\
0\\
\end{array}
\right),
\]
respectively. The dipoles are also members of the family \eqref{E:family} and are given by
\[
\omega_{dipole}(x,y,t) = e^{-\frac{\nu}{\delta^2} t} \sin (x/\delta) + e^{-\nu t}\sin y,
\]
or similarly with sine replaced by cosine, with velocity field
\[
\uu_{dipole}(x,y,t) =  \left(
\begin{array}{c}
e^{-\nu t}\cos y\\
-\delta e^{-\frac{\nu}{\delta^2}t}\cos(x/\delta)\\
\end{array}
\right).
\]

\begin{figure}[H]
\centering
\begin{subfigure}[b]{0.32\textwidth}
    \centering
  \includegraphics[width=\linewidth]{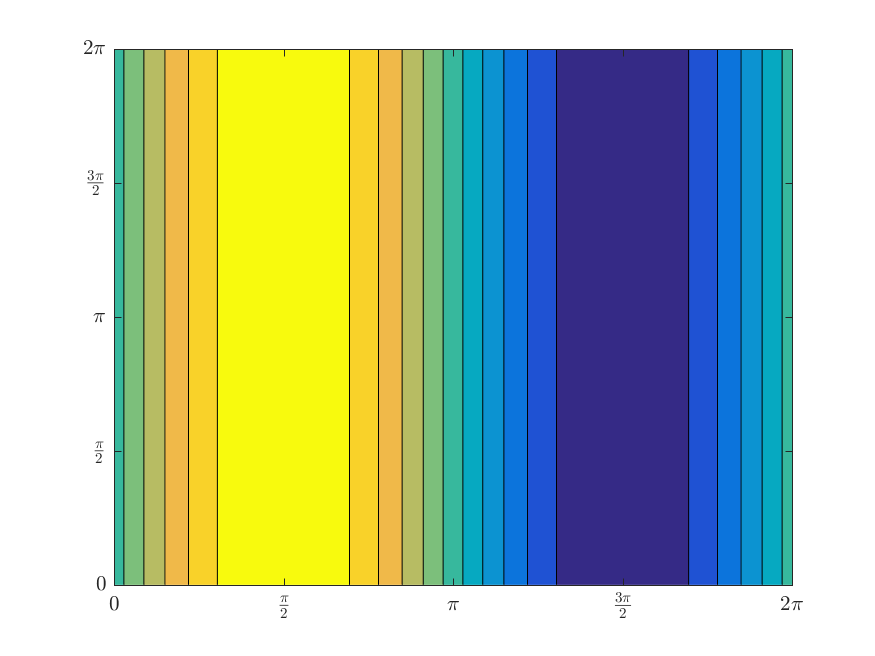}
  \caption{x-bar: $\omega_{xbar}=\sin(x)$}\label{fig:xbar contour}
\end{subfigure}\hfill
\begin{subfigure}[b]{0.32\textwidth}
\centering
  \includegraphics[width=\linewidth]{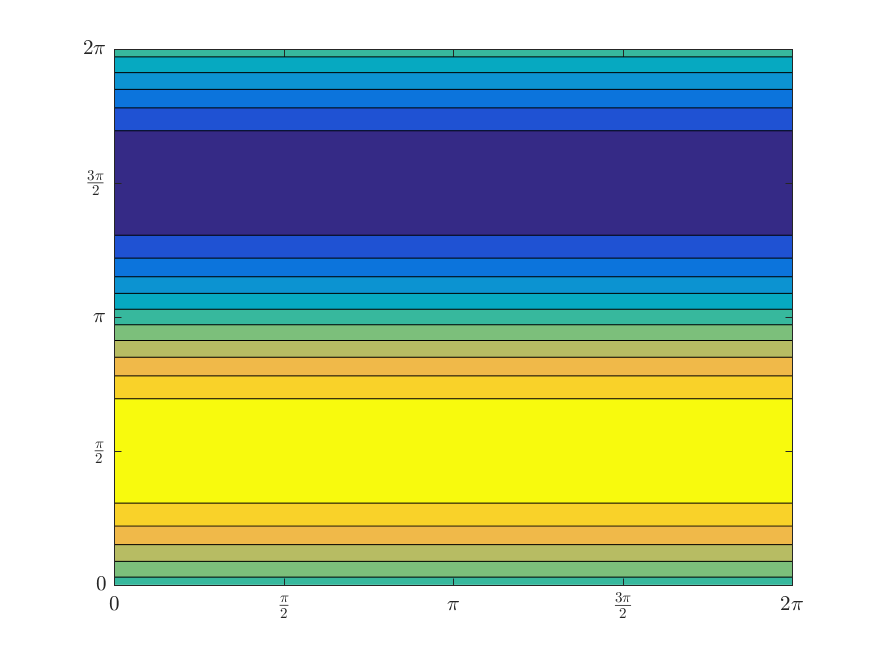}
  \caption{y-bar: $\omega_{ybar}=\sin(y)$ }\label{fig:ybar contour}
\end{subfigure}\hfill
\begin{subfigure}[b]{0.32\textwidth}
\centering
  \includegraphics[width=\linewidth]{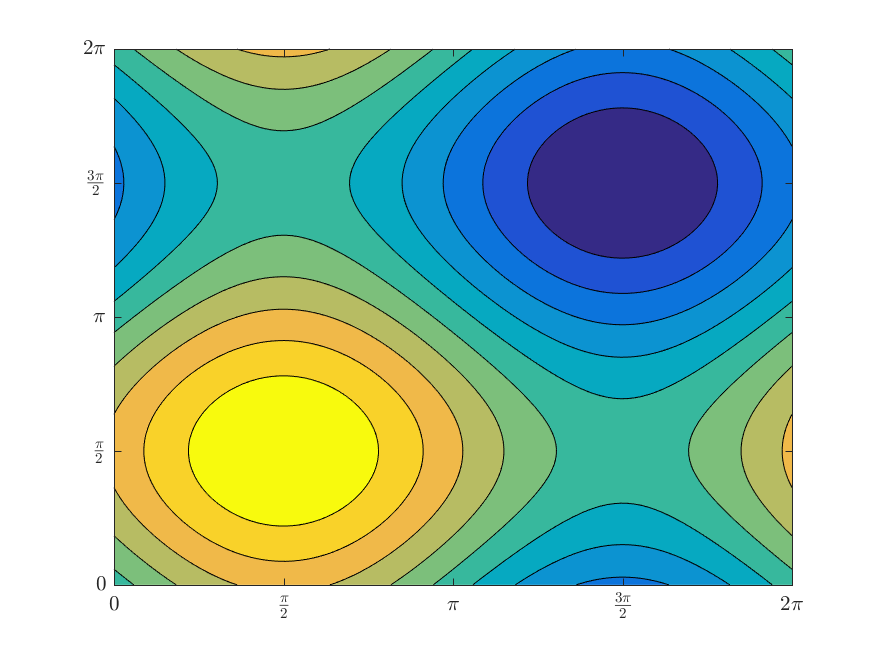}
  \caption{Dipole: $\omega_{dipole}=\sin(x)+\sin(y)$}\label{fig:dipole contour}
\end{subfigure}
\caption{Contour plots of the three quasi-stationary states on the symmetric torus}\label{fig:contour plots}
\end{figure}

For illustration, contour plots for the bar and dipole states for fixed $t=0$ on the symmetric torus ($\delta=1$) are shown in Figure \ref{fig:contour plots}.

When $\nu = 0$, equation \eqref{E:2dns} becomes the Euler equation. It is reasonable to expect that stationary solutions of the Euler equation could play an critical role in the evolution of the Navier-Stokes equation for $0 < \nu \ll 1$. However, there are infinitely many stationary solutions, including the bars and dipoles, and so it is not immediately clear how to determine which would be most important. In \cite{Yin}, entropy arguments and extensive numerical studies were conducted in the case $\delta = 1$ and suggested that the bars and dipoles should be the two most important stationary solutions of the Euler equations. Although both states were observed after initial transient periods in the evolution of the Navier-Stokes equation, interestingly the dipole seemed to emerge for a large class of initial data, whereas the bar states only emerged for a special class of initial data. Subsequent work, again for the deterministic system, showed that indeed the bar states attract nearby solutions at a rate much faster than the background global decay rate, confirming their importance as quasi-stationary states. Results in the case $\delta = 1$ can be found in \cite{BeckWayne13, IbrahimMaekawaMasmoudi17} and results for more general values of $\delta$ are in \cite{LinXu, weiz}.
The stochastic system \eqref{E:2dvortS} was numerically analyzed in \cite{BouchetSimonnet09} where, after an initial transient period, metastable switching between the bars and dipoles was seen, with the dipole being dominant for $\delta = 1$ and the bar states being dominant for $\delta \neq 1$.

In this paper we develop a low dimensional model that captures how the dominant quasi-stationary state in the stochastically forced Navier-Stokes equation is selected by the aspect ratio of the spatial domain, $\delta$. Among the existing results, those that most greatly motivate this paper can be found in \cite{BeckCooperSpiliopoulos} and \cite{BouchetSimonnet09}. The results of the latter paper \cite{BouchetSimonnet09}, briefly described above, to our knowledge were the first to suggest that $\delta$ could provide such a selection mechanism. The former paper \cite{BeckCooperSpiliopoulos} was our previous work focusing on the deterministic vorticity equation, \eqref{E:2dvort}, in which we derived a finite-dimensional model that captured the selection mechanism via the parameter $\delta$. We now seek to use that same finite-dimensional model, but with the addition of noise, to numerically investigate the selection mechanism for the stochastic equation \eqref{E:2dvortS}.  Indeed, one can see from the sample paths of Figure \ref{fig:sample_path}, that individual sample paths exhibit transitions between $x$-bar and $y$-bar states, as it has also been observed in \cite{BouchetSimonnet09}.

The rest of the paper is organized as follows. In \S \ref{S:FourierSpaceRep} we review the finite-dimensional model originally derived in \cite{BeckCooperSpiliopoulos} and the theoretical results of that work regarding the selection mechanism in the deterministic setting. In addition, we also add noise to that model to obtain the stochastic differential equation (SDE) model that is the focus of this current work.
In \S\ref{S:sims}, to determine the validity of the SDE model, we compare statistics related to a direct simulation of the stochastic vorticity equation \eqref{E:2dvortS} with those of the SDE. We demonstrate numerically that the statistics of the two equations agree in all cases, $\delta>1$, $\delta<1$ and $\delta=1$. In particular, solutions to both systems evolve towards an $x$-bar, $y$-bar, and dipole in the three respective cases. In \S \ref{S:PDE}, we further examine the SDE model by viewing it as a perturbation in the limit as $|\delta^2-1|$ and $\nu$ converge to zero. We show that, after appropriate time-space rescalings, the system can be viewed as a slow-fast system and classical averaging and homogenization techniques apply. Via the backward Kolmogorov equation, a system of PDEs that governs the leading order dynamics of a key order parameter, $\mathbb{E}[Z_{red}(t)]$, defined in \eqref{E:op}, is derived. This gives us an additional formal approximation to the expected value of the order parameter, which we can use to show the selection of the quasi-stationary state. Numerically solving the PDEs allows us to approximate the evolution of $\mathbb{E}[Z_{red}(t)]$ for values of $\delta$ close to 1, at least on some initial finite interval of time. Conclusions and future directions then follow in \S \ref{S:Conclusion}.


\section{Fourier space representation and model reduction}\label{S:FourierSpaceRep}

Due to the form of the family of solutions \eqref{E:family}, it is most convenient to express the stochastic vorticity equation in Fourier space.  Hence, letting
\[
\omega(x,y) = \sum_{\vec{k} \neq (0,0)} \hat{\omega}_{\vec{k}} e^{\rmi (k_1x/\delta + k_2y)}, \qquad \hat{\omega}_{\vec{k}} = \frac{1}{4\pi^2\delta}\int_{D_\delta} \omega(x,y) e^{-\rmi(k_1x/\delta + k_2y)} \rmd x \rmd y,
\]
we obtain, for $\vec{j}$, $\vec{k}$ and $\vec{l}\neq(0,0)$, the following system of infinitely many coupled SDEs,
\begin{equation}
\begin{aligned}
\dot{\hat{\omega}}_{\vec{k}} &= -\frac{\nu}{\delta^2} |\vec{k}|^2_\delta \hat{\omega}_{\vec{k}} - \delta \sum_{\vec{l}} \frac{\langle \vec{k}^\perp, \vec{l}\rangle}{|\vec{l}|_\delta^2} \hat{\omega}_{\vec{k}-\vec{l}}\hat{\omega}_{\vec{l}} +\sqrt{2\nu}\sigma_{\vec{k}}\dot{\beta}_{\vec{k}}\\
&= -\frac{\nu}{\delta^2} |\vec{k}|^2_\delta \hat{\omega}_{\vec{k}} - \frac{\delta}{2} \sum_{\vec{j}+\vec{l}=\vec{k}} \langle \vec{j}^\perp, \vec{l}\rangle \left( \frac{1}{|\vec{l}|_\delta^2} -  \frac{1}{|\vec{j}|_\delta^2} \right) \hat{\omega}_{\vec{j}}\hat{\omega}_{\vec{l}}+\sqrt{2\nu}\sigma_{\vec{k}}\dot{\beta}_{\vec{k}}, \label{E:vortS}
\end{aligned}
\end{equation}
where
\begin{equation}\label{E:norm}
|\vec{k}|_\delta^2 = k_1^2 + \delta^2 k_2^2, \qquad \vec{k}^\perp = (k_2, -k_1).
\end{equation}

Viewing the system in Fourier space allows us to use the relative energy in certain modes to measure the proximity of solutions to an $x$-bar, $y$-bar, or dipole state.  The $x$-bar states, $e^{-\frac{\nu}{\delta^2} t}\cos(x/\delta)$ and $e^{-\frac{\nu}{\delta^2} t}\sin(x/\delta)$, correspond to solutions with energy only in the $\vec{k}=(\pm1,0)$ modes and the $y$-bar states, $e^{-\nu t}\cos(y)$ and $e^{-\nu t}\sin(y)$, correspond to solutions with energy only in the $\vec{k}=(0,\pm1)$ modes. Solutions with energy in both the $\vec{k}=(\pm1,0)$ and $\vec{k}=(0,\pm1)$ modes correspond to the dipole state. These four modes are the lowest in the system and will be referred to as the ``low modes". They correspond to modes with the lowest value of $|\vec{k}|_\delta$ defined by \eqref{E:norm}.
Any mode $\hat{\omega}_{\vec{k}}$ with $|\vec{k}|>\mbox{max}\{1,\delta^2\}$ will from here on be referred to as a ``high mode".

To measure the relative energy in the low modes, we define the stochastic order parameter,
\begin{equation}\label{E:op_stochastic}
    Z_{vort}(t):=\frac{|\hat{\omega}_{(1,0)}(t)|^2}{|\hat{\omega}_{(1,0)}(t)|^2+|\hat{\omega}_{(0,1)}(t)|^2},
\end{equation}
where $\hat{\omega}_{(1,0)}$ and $\hat{\omega}_{(0,1)}$ solve \eqref{E:vortS}. Due to the condition,  $\hat{\omega}_{(k_1,k_2)}=\bar{\hat{\omega}}_{(-k_1,-k_2)}$, the relative energy in all of the low modes can be captured by $Z_{vort}(t)$. The value of $Z_{vort}(t)$, bounded between 0 and 1, corresponds to the proximity of the solution to an $x$-bar, $y$-bar or dipole state. If the dynamics drive $Z_{vort}(t)$ to increase to 1, there is more energy in $\hat\omega_{(1,0)}$ relative to $\hat\omega_{(0,1)}$, indicating the system is in an $x$-bar state. Conversely if $Z_{vort}(t)$ falls toward 0, the system would be observed to be in a $y$-bar state. If $Z_{vort}(t)$ instead stays near 1/2, the system is in a dipole state with relative energy in the low modes comparable in magnitude.

The finite dimensional system that we will use to model \eqref{E:vortS} will be defined in terms of the lowest eight Fourier modes, which for notational convenience we denote as
\begin{eqnarray}
\omega_1 &:=& \hat\omega_{(1,0)},\quad \omega_2 := \hat\omega_{(-1, 0)},\quad \omega_3 := \hat\omega_{(0, 1)}, \quad \omega_4 := \hat\omega_{(0, -1)}, \nonumber \\
\omega_5 &:=& \hat\omega_{(1,1)},\quad \omega_6 := \hat\omega_{(-1, 1)}, \quad \omega_7 := \hat\omega_{(1, -1)}, \quad  \omega_8 := \hat\omega_{(-1, -1)}. \label{E:fourier}
\end{eqnarray}

The variables $\omega_{1, 2, 3, 4}$ correspond to the low modes, while $\omega_{5, 6, 7, 8}$ represent the role of all the high modes. Since the solution $\omega(x,y)$ of \eqref{E:vortS} is real valued, the following complex conjugacy relationship must still hold,
\begin{equation}\label{symmetry}
\omega_1 = \bar \omega_2, \quad \omega_3 = \bar \omega_4, \quad \omega_5 = \bar \omega_8, \quad \omega_7 = \bar \omega_8.
\end{equation}

Thus the reduced model will be an eight dimensional approximation to the dynamics of \eqref{E:vortS}. To derive the model, we apply a center manifold reduction to \eqref{E:vortS} with $\sigma_{\vec{k}} = 0$ for all $\vec{k}$ to obtain an eight-dimensional deterministic ODE, which is the model studied in \cite{BeckCooperSpiliopoulos}, and then add noise back to that system to obtain the final eight-dimensional SDE model we study here.

To carry out the center manifold reduction onto the lowest eight modes, assume for $\hat{\omega}_{\vec{k}}$ with $\vec{k} \notin \{(\pm1, 0), (0, \pm1), (\pm1, \pm1)\} =: \mathcal{K}_0$, that there exists a smooth function $H(\omega_1, \dots, \omega_8; \vec{k})$ such that the eight-dimensional manifold defined by
\[
\mathcal{M} = \{ \hat{\omega}: \hat{\omega}_{\vec{k}} = H(\omega_1, \dots, \omega_8; \vec{k}), \quad \vec{k} \notin \mathcal{K}_0\}
\]
is invariant for the deterministic dynamics of \eqref{E:vortS} with $\sigma_{\vec{k}}=0$ for all $\vec{k}$. We refer to this as a center manifold because it is defined in terms of the lowest eight modes, which have the weakest linear decay rates. Based on this assumption, one can then in principle compute the coefficients of the
Taylor expansion of $H(\cdot, \vec{k})$ to any order for each $\vec{k}$ by taking the derivative of each of the low modes in two ways (via the function $H$ and \eqref{E:vortS} with $\sigma_{\vec{k}=0}$) and equating coefficients. See \cite{BeckCooperSpiliopoulos} for the details of the derivation.

The reduction is local and will only be valid in a size $\mathcal{O}(\nu)$ neighborhood of zero due to the small spectral gaps for the operator $\nu\Delta$. Additionally, while the existence of a finite dimensional (inertial) model of the system \eqref{E:vortS} that describes the global dynamics cannot be expected \cite{Zelik14}, the model still provides meaningful insight into the role $\delta$ plays in selecting the dominant quasi-stationary state for small initial conditions. For additional examples in which similar reductions of the Navier-Stokes equation to a finite dimensional model have been used to understand global dynamics see \cite{EMattingly01,MattinglyPardoux14}.

Adding independent (real) Brownian motions $W_{1,3,5,7}$ to each equation of the resulting ODE model leads to our final SDE model
\begin{equation}
    \begin{aligned}
        \dot \omega_1 &= - \frac{\nu}{\delta^2} \omega_1 + \frac{1}{\delta(1+\delta^2)}[\omega_3\omega_7 - \bar{\omega}_3 \omega_5] + \frac{3\delta^6}{2\nu(4+\delta^2)(1+\delta^2)^2}\omega_1(|\omega_5|^2 + |\omega_7|^2) +\sqrt{2\nu}\sigma_1\dot{W}_1 \\
        \dot \omega_3 &= - \nu \omega_3 + \frac{\delta^3}{(1+\delta^2)}[\bar{\omega}_1\omega_5 - \omega_1 \bar{\omega}_7] + \frac{3\delta^2}{2\nu(1+4\delta^2)(1+\delta^2)^2}\omega_3(|\omega_5|^2 + |\omega_7|^2) +\sqrt{2\nu}\sigma_3\dot{W}_3\label{E:deltanot1}  \\
        \dot \omega_5 &= - \nu \frac{1+\delta^2}{\delta^2} \omega_5 -\frac{\delta^2-1}{\delta}\omega_1 \omega_3   -\frac{\delta^6(3+\delta^2)}{2\nu(4+\delta^2)(1+\delta^2)} \omega_5 |\omega_1|^2-\frac{1+3\delta^2}{2\nu\delta^2(1+4\delta^2)(1+\delta^2)}\omega_5 |\omega_3|^2  +\sqrt{2\nu}\sigma_5\dot{W}_5 \\
        \dot \omega_7 &= - \nu \frac{1+\delta^2}{\delta^2} \omega_7 +\frac{\delta^2-1}{\delta}\omega_1 \bar{\omega}_3 -\frac{\delta^6(3+\delta^2)}{2\nu(4+\delta^2)(1+\delta^2)} \omega_7 |\omega_1|^2 - \frac{1+3\delta^2}{2\nu\delta^2(1+4\delta^2)(1+\delta^2)}\omega_7 |\omega_3|^2 +\sqrt{2\nu}\sigma_7\dot{W}_7.
    \end{aligned}
\end{equation}

Note that \eqref{E:deltanot1} with $\sigma_{1,3,5,7}=0$ corresponds to the ODE model derived in \cite{BeckCooperSpiliopoulos}. To compare the dynamics of this model to that of $Z_{vort}(t)$, defined in \eqref{E:op_stochastic}, we define the analogous order parameter for the SDE model,
\begin{equation}\label{E:op}
    Z_{red}(t):=\frac{|\omega_1(t)|^2}{|\omega_1(t)|^2+|\omega_3(t)|^2},
\end{equation}
which again is used to determine towards which quasi-stationary state the system trends. Here, $\omega_1(t)$ and $\omega_3(t)$ are solutions to the reduced system \eqref{E:deltanot1}. The Monte Carlo simulation of the reduced model finds that the dominant quasi-stationary state depends on the aspect ratio of $D_\delta$ in the same way as the deterministic model, studied in detail in \cite{BeckCooperSpiliopoulos}. The main result there, which describes the selection of quasi-stationary states in \eqref{E:deltanot1} with $\sigma_{1,3,5,7}=0$, can be described by the following theorem.

\begin{Theorem}\cite[Theorem 3.4]{BeckCooperSpiliopoulos}\label{T:det}
For $\delta\in\left(\sqrt{\frac{2}{3}},\sqrt{\frac{3}{2}}\right)$, under the dynamics of \eqref{E:deltanot1} with $\sigma_{1,3,5,7}=0$, if $\delta>1$, then $Z_{red}(t)\rightarrow1$, indicating evolution to an $x$-bar state. Conversely if $\delta<1$, then $Z_{red}(t)\rightarrow0$, indicating evolution to a $y$-bar state. For $\delta=1$, there exists a one-dimensional center manifold of fixed points in the phase space that determines the asymptotic limit of $Z_{red}(t)$. This center manifold is foliated with co-dimension one stable manifolds in which solutions converge to the corresponding fixed point. Exactly one of these manifolds corresponds to each of the limits $Z_{red}(t)\rightarrow1$ and $Z_{red}(t)\rightarrow0$. Thus, generic initial conditions are seen to evolve to the dipole state.
\end{Theorem}

\begin{Remark}
The order parameter considered in \cite{BeckCooperSpiliopoulos} was instead the ratio $R(t)=|\omega_1(t)|^2/|\omega_3(t)|^2$. Theorem \ref{T:det} frames the result in terms of the order parameter $Z_{red}(t)$. The choice to now consider $Z_{red}(t)$ is for convenience with regards to numerical simulation due to its being bounded between 0 and 1.
\end{Remark}

\begin{Remark}\label{R:RealSystem}
A straightforward computation shows that, for any $\delta$, the set $\{\Im(\omega_1)=\Im(\omega_3)=\Im(\omega_5)=\Im(\omega_7)=0\}$ is invariant under the dynamics of \eqref{E:deltanot1} with $\sigma_{1,3,5,7}=0$. Since the real subsystem is invariant in the deterministic setting, we simulate the reduced model where the modes, $\omega_{1,3,5,7}$, as well as the Wiener processes, $W_{1,3,5,7}$, are all real valued.
\end{Remark}

\section{Numerical simulation of the vorticity equation and reduced model}\label{S:sims}

This section provides simulations of the vorticity equation \eqref{E:vortS} and of the reduced model \eqref{E:deltanot1}. Via Monte Carlo simulation, the average evolution of the order parameters $Z_{vort}(t)$ and $Z_{red}(t)$ will be plotted for several values of $\delta$ near 1. It will be seen that the reduced model captures the selection of the quasi-stationary states via the parameter $\delta$. In particular, in both models, for a particular value of $\delta\approx1$, the system's selection of its dominant quasi-stationary state is consistent with the motivating results, given by Theorem \ref{T:det}. In particular, the system selects, as the dominant quasi-stationary state, a dipole for $\delta=1$, an $x$-bar for $\delta>1$, and a $y$-bar for $\delta<1$.

The simulation of \eqref{E:vortS} is done via a spectral method which includes Fourier modes $\hat{\omega}_{\vec{k}}$ with $\vec{k}\in\mathcal{K}:=\{\vec{k}=(k_1,k_2)\in \mathbb{Z}^2:   0\leq|k_1|,|k_2|\leq64 \mbox{ and } (k_1,k_2)\neq (0,0)  \}$; see \cite{Lord}. A condition of exponential decay is imposed on the noise coefficients $\sigma_{\vec{k}}$ seen in \eqref{E:noise},
\begin{equation}\label{E:noise-decay}
|\sigma_{\vec{k}}|\leq e^{-\alpha_0|\vec{k}|^2}.
\end{equation}

Similar to \cite{BouchetSimonnet09}, simulations are conducted with
$\sum_{\{\vec{k}\in\mathcal{K}\}}e^{-\alpha_0|\vec{K}|^2}=1$. For our
set $\mathcal{K}$, this means $\alpha_0\approx0.349$. Time was finely
discretized and a tamed semi-implicit Euler-Maruyama method was implemented to simulate the stochastically forced reduced system.

We verify the selection of the dominant quasi-stationary state using Monte Carlo simulation where the average path over $N$ trials is plotted. Individual runs will be denoted by $Z^i_{vort}(t)$ and $Z^i_{red}(t)$, for $i=1,\dots N$, with corresponding averages given by
\[
\bar{Z}_{vort}(t)=\frac{1}{N}\sum_{i=1}^{N}Z^{i}_{vort}(t), \qquad \bar{Z}_{red}(t)=\frac{1}{N}\sum_{i=1}^{N}Z^{i}_{red}(t).
\]

Similarly, we define the empirical variances to be
\[
V_{vort}(t)=\frac{1}{N-1}\sum_{i=1}^{N}(Z^{i}_{vort}(t)-\bar{Z}_{vort}(t))^2, \qquad V_{red}(t)=\frac{1}{N-1}\sum_{i=1}^{N}(Z^{i}_{red}(t)-\bar{Z}_{red}(t))^2.
\]

It will also be useful to plot the time averages of these Monte Carlo averages. To produce a meaningful average we introduce a ``burn-in time", $t_{burn}$, and ignore the initial period during which $\bar{Z}_{vort}(t)$ and $\bar{Z}_{red}(t)$ have not yet stabilized. Define this time average for any function $f(t)$ defined on $t_{burn}\leq t\leq T$ to be
\[
A(f,t_{burn}):=\frac{1}{T-t_{burn}}\int_{t_{burn}}^{T} f(t) \mbox{ } \rmd t.
\]

\subsection{Vorticity Equation}

\begin{figure}[H]
        \centering
        \begin{subfigure}[b]{0.475\textwidth}
            \centering
            \includegraphics[width=\textwidth]{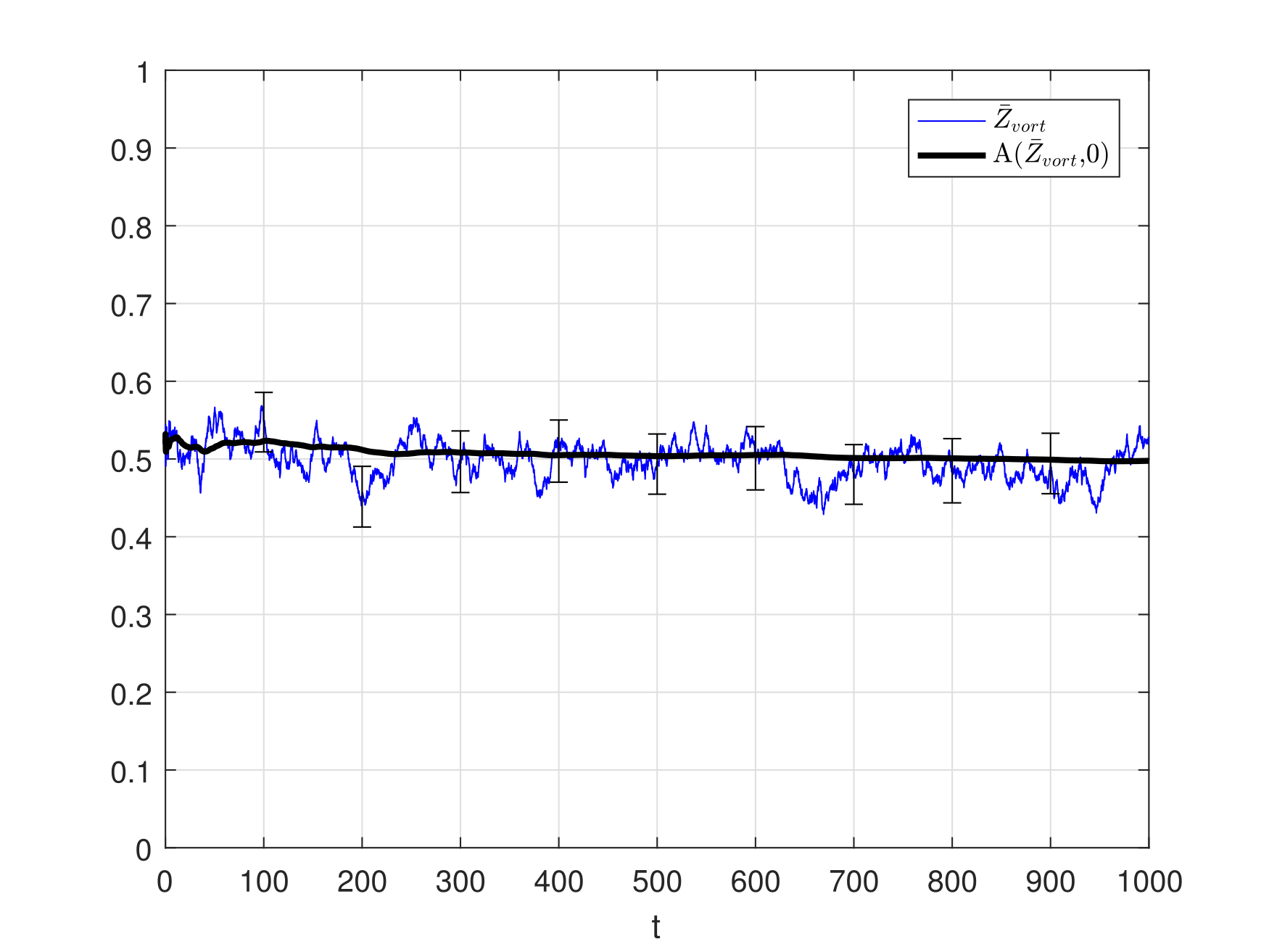}
            \caption{$\bar{Z}_{vort}(t)$ with 95\% confidence interval.}
            \label{fig:delta 1 vort}
        \end{subfigure}
        \hfill
        \begin{subfigure}[b]{0.475\textwidth}
            \centering
            \includegraphics[width=\textwidth]{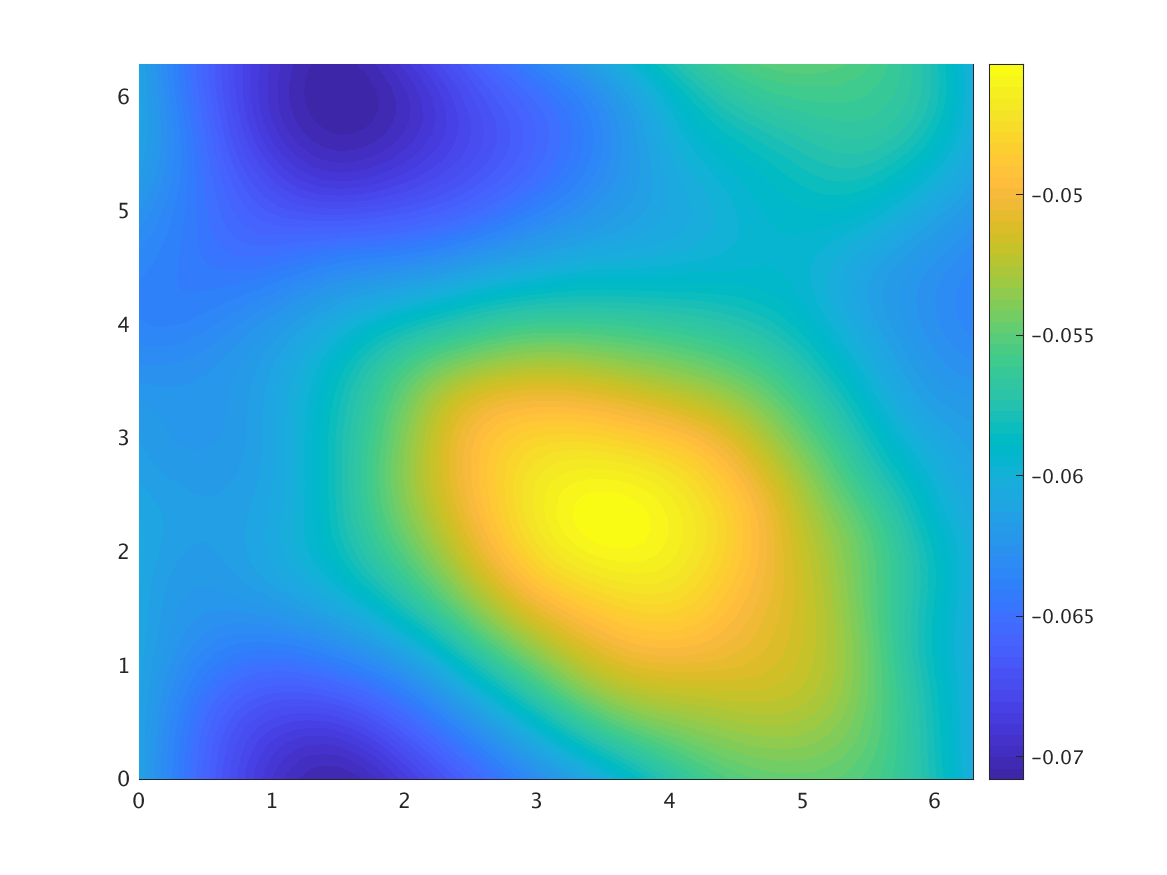}
            \caption{Average contour plot of vorticity.}
            \label{fig:delta 1 vort contour}
        \end{subfigure}
                \caption{Vorticity aligns on average as a dipole for $\delta=1$.}\label{fig:delta 1 CI}
\end{figure}

\begin{figure}[H]
        \centering
        \begin{subfigure}[b]{0.475\textwidth}
            \centering
            \includegraphics[width=\textwidth]{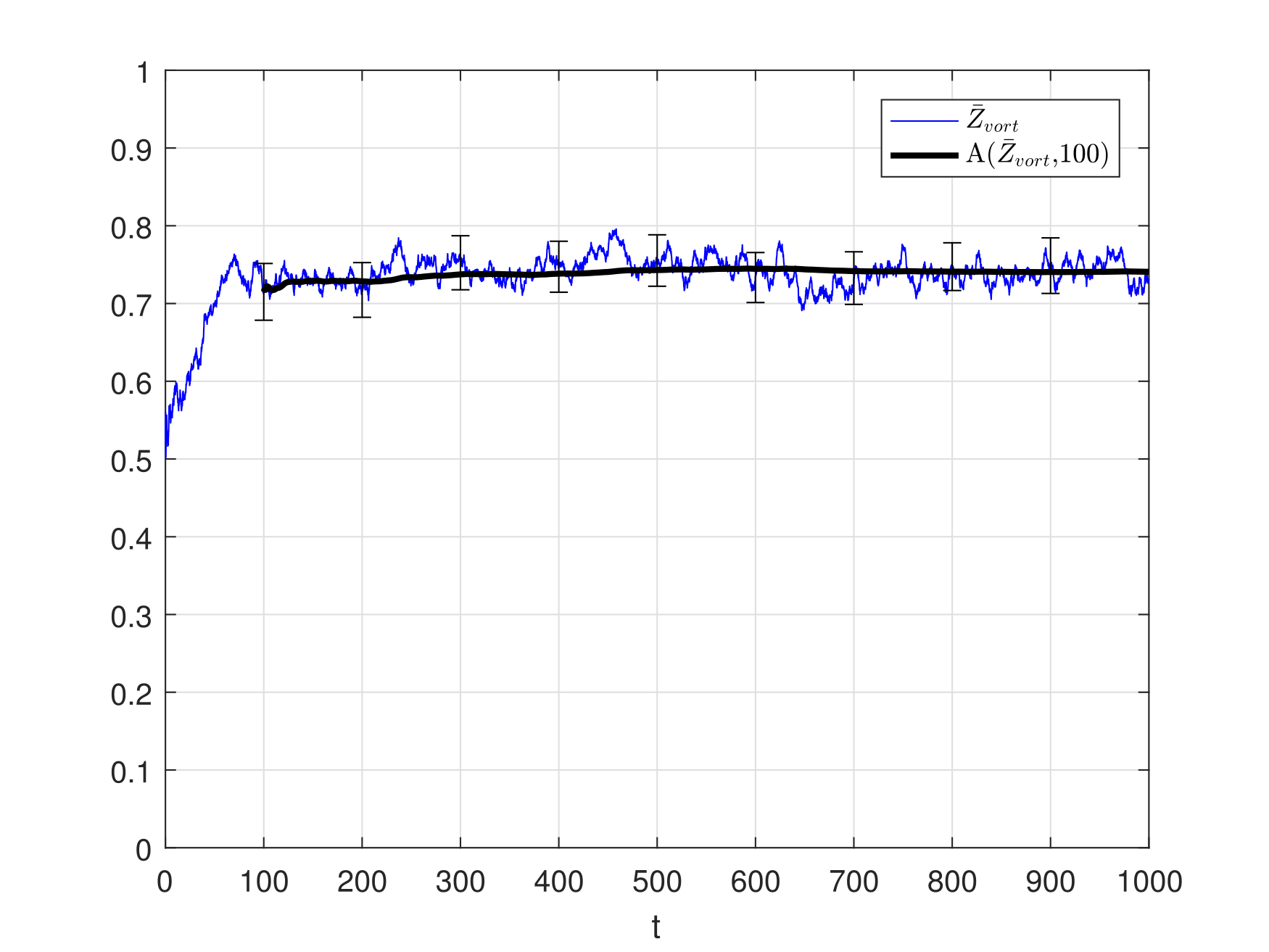}
            \caption{$\bar{Z}_{vort}(t)$ with 95\% confidence interval.}
            \label{fig:xbar full}
        \end{subfigure}
        \hfill
        \begin{subfigure}[b]{0.475\textwidth}
            \centering
            \includegraphics[width=\textwidth]{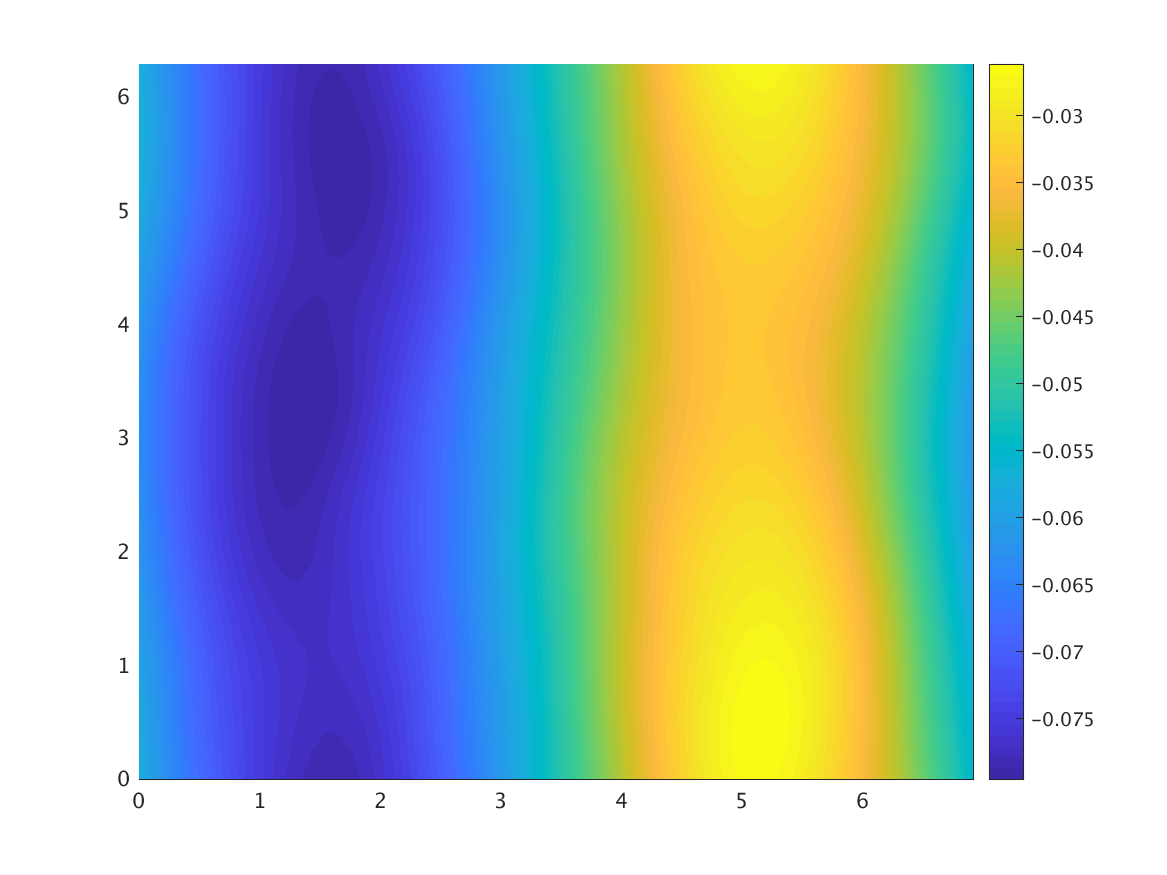}
            \caption{Contour plot of vorticity.}
            \label{fig:xbar full contour}
        \end{subfigure}
                \caption{Vorticity aligns on average as an $x$-bar for $\delta=1.1$.}\label{fig:delta 1.1 CI}
\end{figure}

\begin{figure}[H]
        \centering
        \begin{subfigure}[b]{0.475\textwidth}
            \centering
            \includegraphics[width=\textwidth]{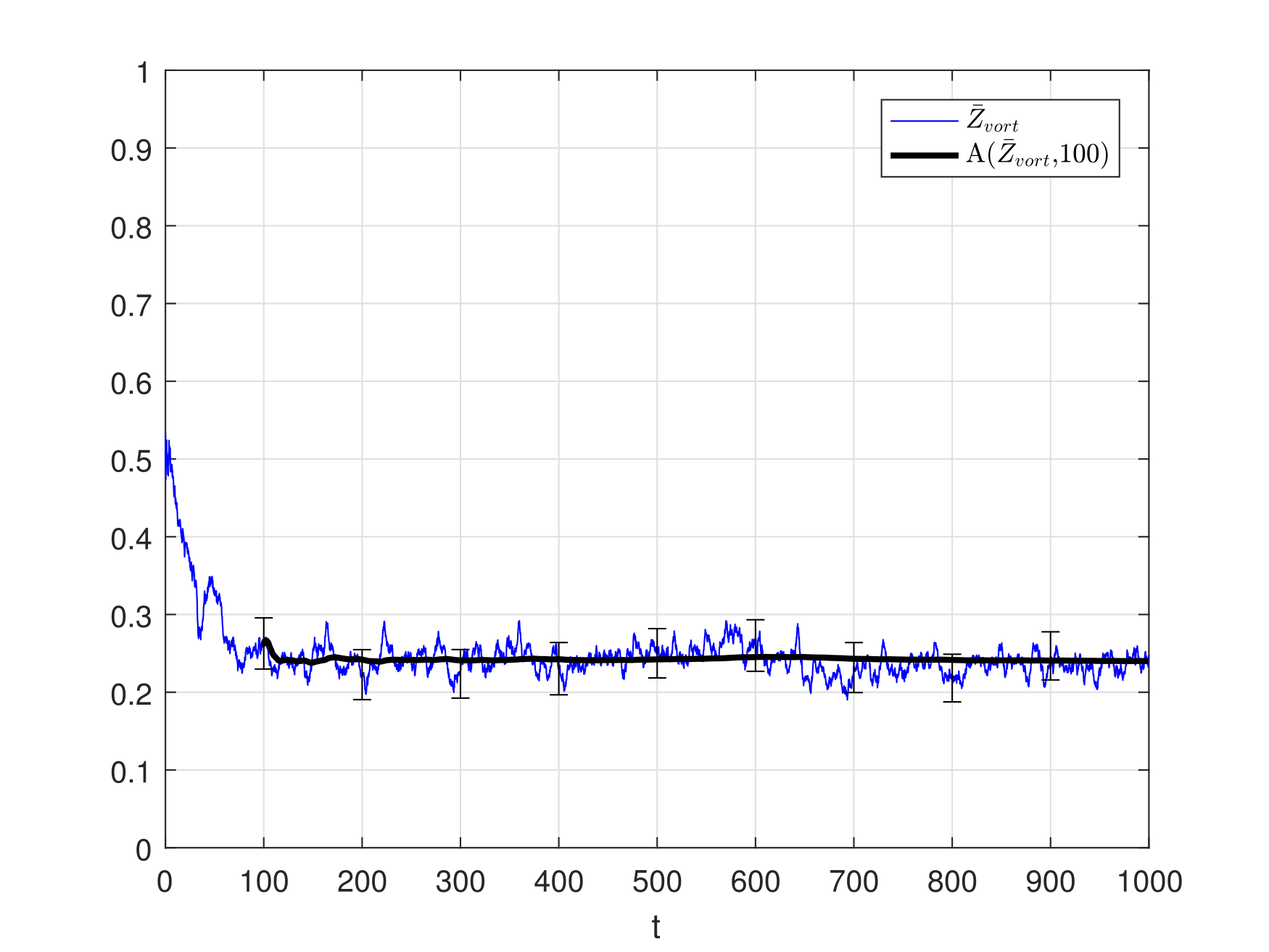}
            \caption{$\bar{Z}_{vort}(t)$ with 95\% confidence interval.}
            \label{fig:ybar full}
        \end{subfigure}
        \hfill
        \begin{subfigure}[b]{0.475\textwidth}
            \centering
            \includegraphics[width=\textwidth]{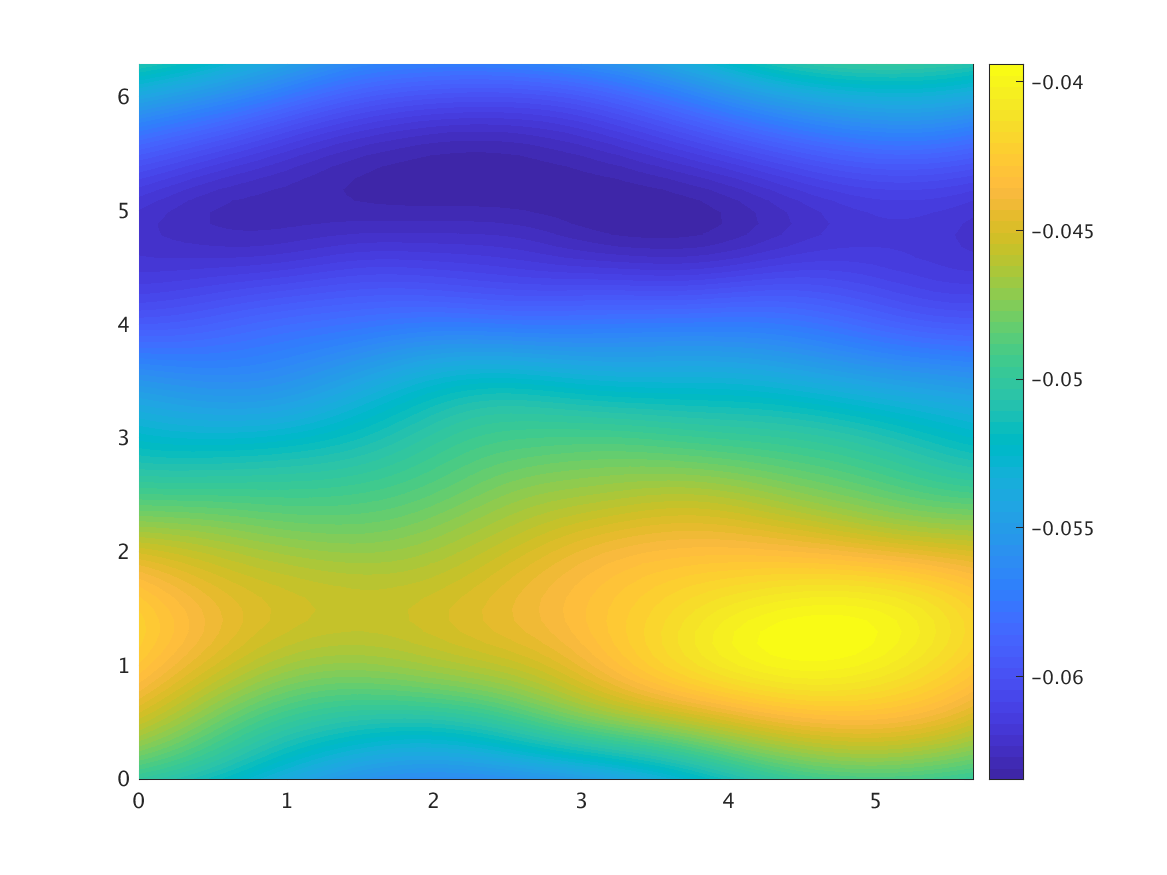}
            \caption{Average contour plot of vorticity.}
            \label{fig:ybar full contour}
        \end{subfigure}
                \caption{Vorticity aligns on average as a $y$-bar for $\delta=0.9$.}\label{fig:delta 0.9 CI}
\end{figure}

Plotted in Figures \ref{fig:delta 1 CI}-\ref{fig:delta 0.9 CI} are $\bar{Z}_{vort}(t)$, the time average $A(\bar{Z}_{vort},t_{burn})$, and the 95\% confidence intervals defined via
\[
CI^{\pm}(t)=\bar{Z}_{vort}(t)\pm 1.96 * \sqrt{\frac{V_{vort}(t)}{N}}.
\]

Also included are average contour plots for the vorticity. We use $N=200$ and for each trial use zero initial conditions and $\nu=0.001$. For $\delta=1$, Figure \ref{fig:delta 1 vort} shows $\bar{Z}_{vort}(t)$ remains near 1/2 for the duration of the simulation. We use a burn-in time of $t_{burn}=0$ when computing the time average since on the symmetric domain it is clear there is no transient initial period. In Figure \ref{fig:delta 1 vort contour}, the average contour plot for each individual trial are themselves averaged over the $N=200$ trials, reflecting a dipole.

The simulations exhibited in Figures \ref{fig:xbar full} and \ref{fig:xbar full contour} show that, for $\delta=1.1$, the order parameter increases initially and the average contour plot looks like that of an $x$-bar state. In Figure \ref{fig:xbar full}, $t_{burn}=100$ is used when computing the time average.

Lastly for $\delta<1$ the simulations exhibited in Figures \ref{fig:ybar full} and \ref{fig:ybar full contour} show that, for $\delta=0.9$,   the order parameter decreases over an initial period of time and the average contour plot looks like that of a $y$-bar state. Here we again set $t_{burn}=100$.

\begin{figure}[H]
        \centering
        \begin{subfigure}[b]{0.32\textwidth}
            \centering
            \includegraphics[width=\textwidth]{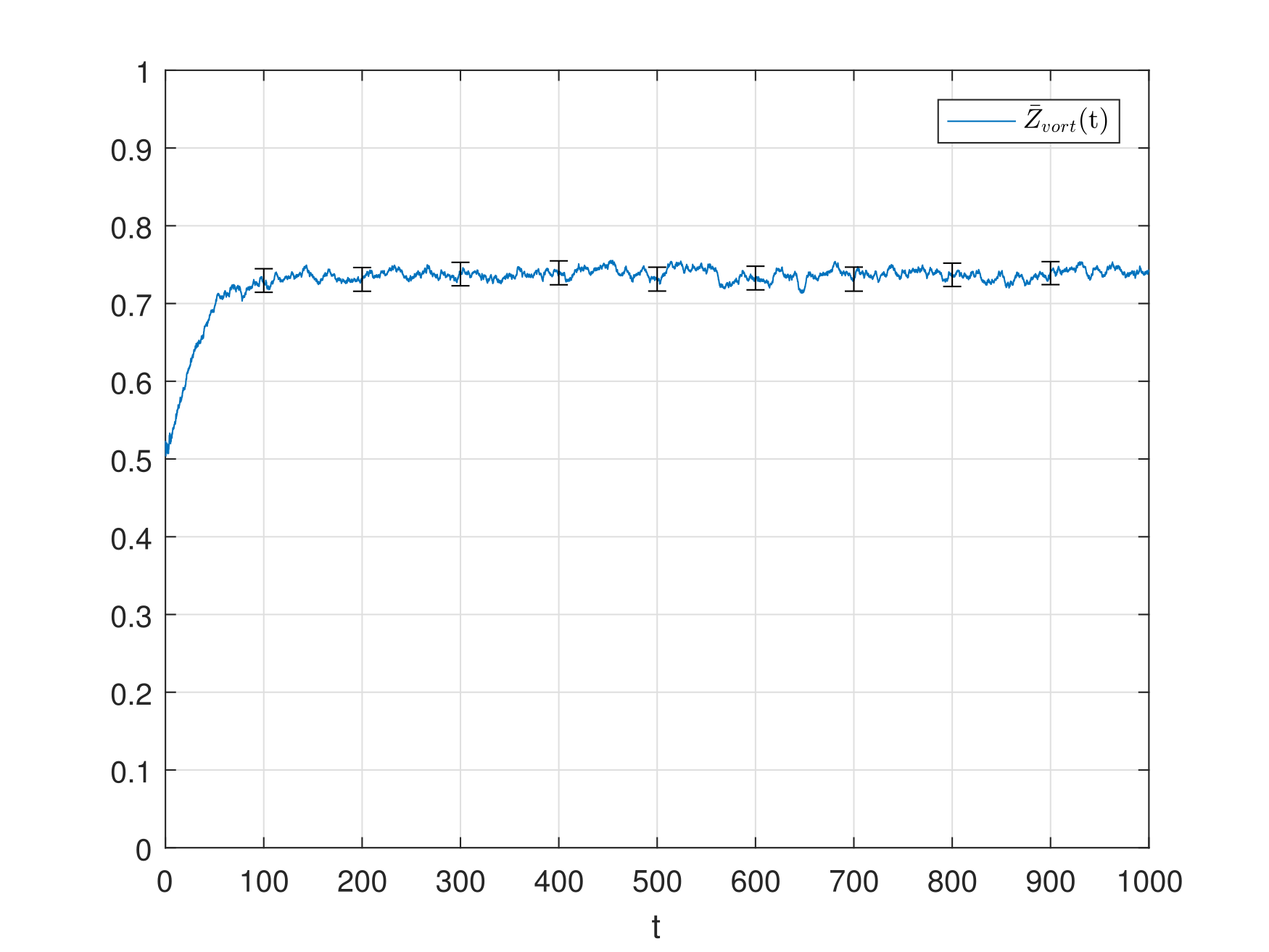}
            \caption{$\delta=1.10$ }
            \label{fig:n1k d110}
        \end{subfigure}
        \begin{subfigure}[b]{0.32\textwidth}
            \includegraphics[width=\textwidth]{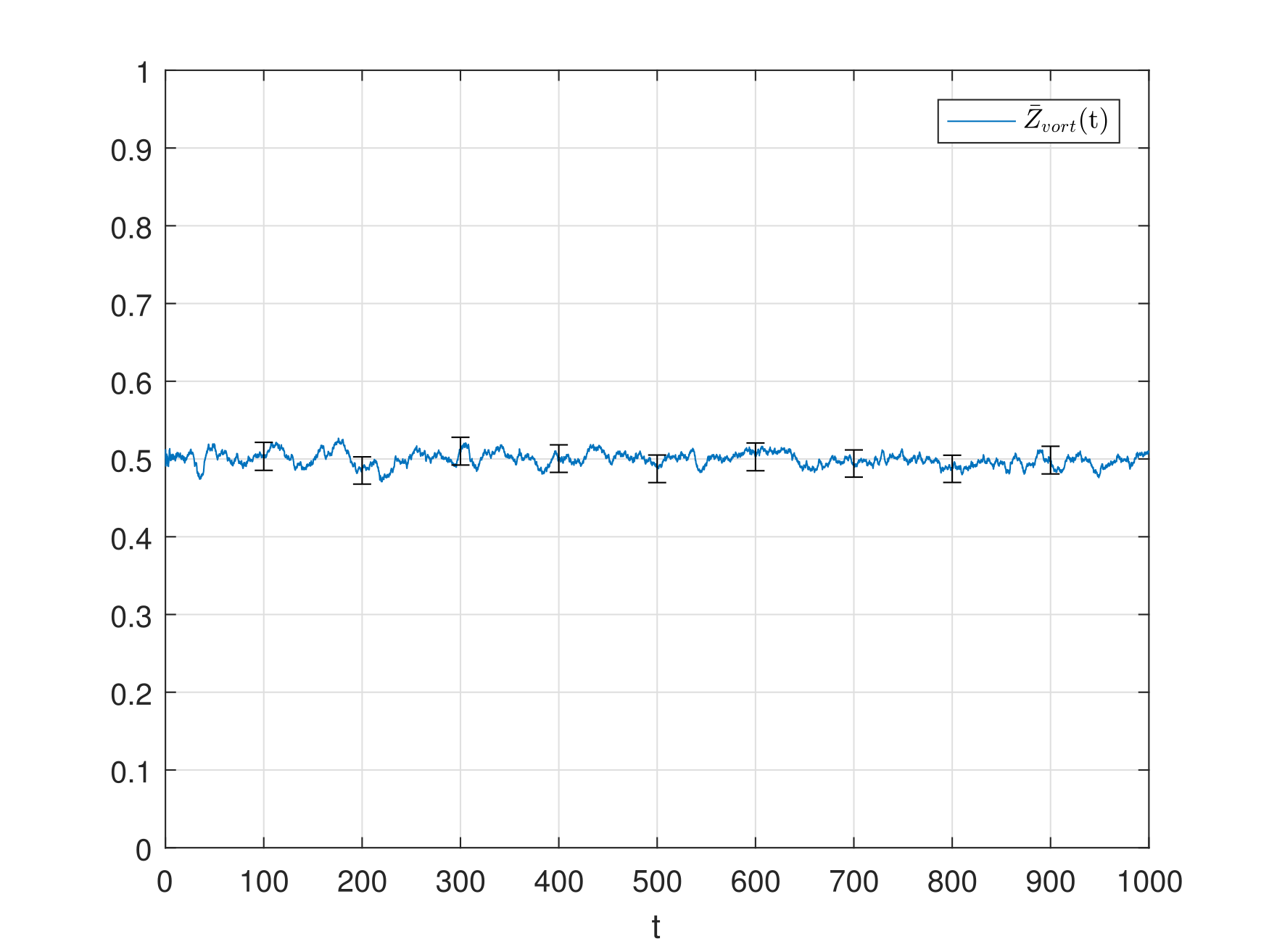}
            \caption{$\delta=1.0$}
            \label{fig:n1k d1}
        \end{subfigure}
          \begin{subfigure}[b]{0.32\textwidth}
            \includegraphics[width=\textwidth]{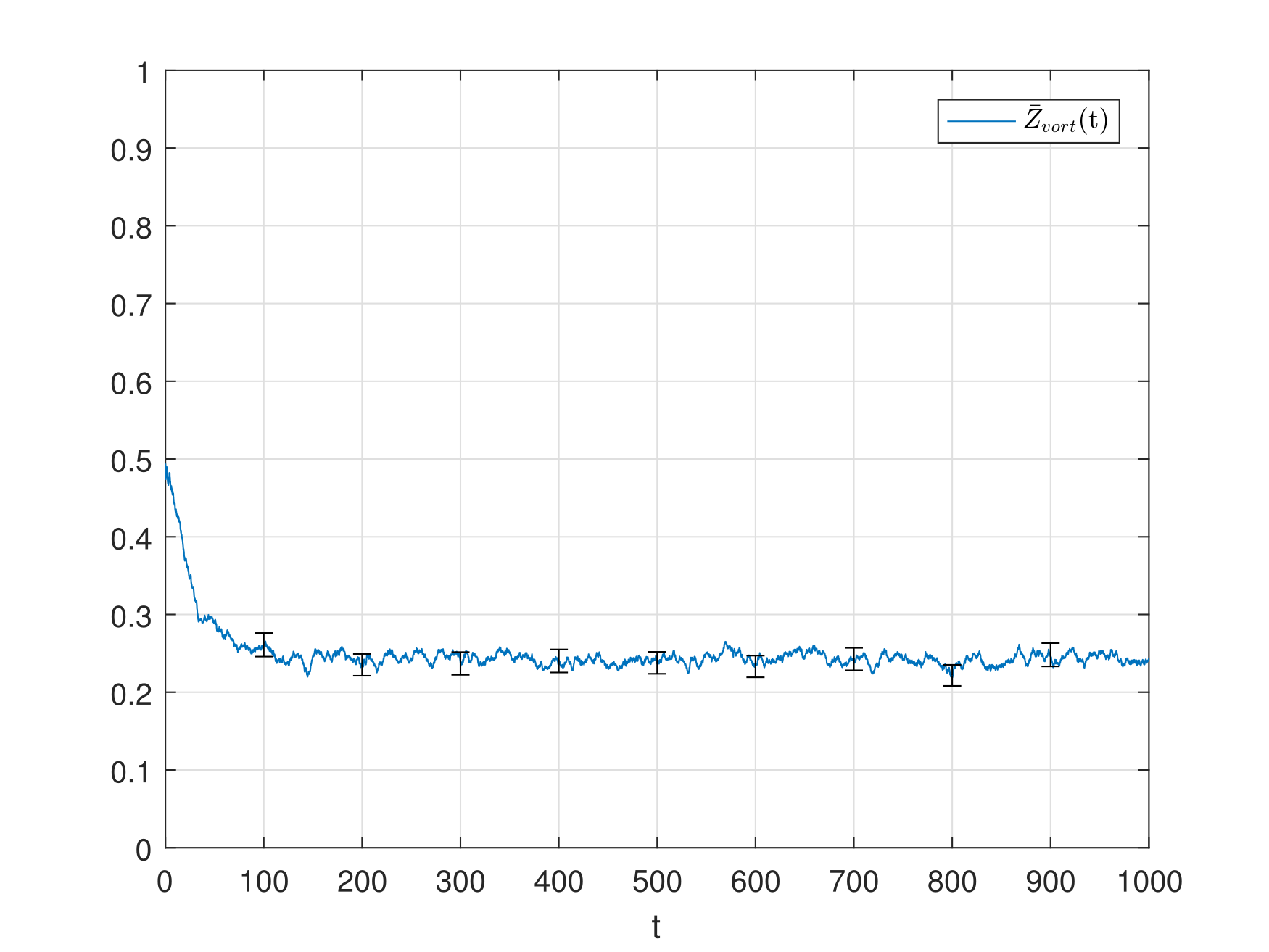}
            \caption{$\delta=0.90$}
            \label{fig:n1k d90}
        \end{subfigure}
                \caption{Plot of  $\bar{Z}_{vort}(t)$ and of 95\% confidence level error bars with $N=1000$ trials and $\nu=0.001$. }
        \label{fig: Zvort n=1000}
\end{figure}

Provided in Figure \ref{fig: Zvort n=1000} are plots of
$\bar{Z}_{vort}(t)$ for $\delta=1.10$, $\delta=1.0$ and $\delta=0.90$
averaged over $N=1000$ trials. This is to show that as the number of
trials increase, the variance is decreasing without changing the mean
behavior.
The variances all remain generally between 0.06-0.08. While the variance does decrease, the limiting value of $\bar{Z}_{vort}(t)$ remains relatively unchanged compared to what is seen when averaging over $N=200$ trials.

\begin{figure}[H]
\begin{subfigure}[b]{0.48\textwidth}
  \includegraphics[width=\linewidth]{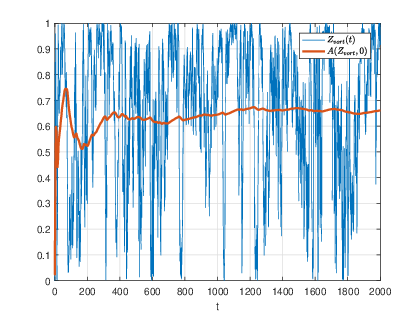}
  \caption{An individual trajectory transitions among quasi-stationary states.}\label{fig:sample_path}
\end{subfigure}
\begin{subfigure}[b]{0.48\textwidth}%
  \includegraphics[width=\linewidth]{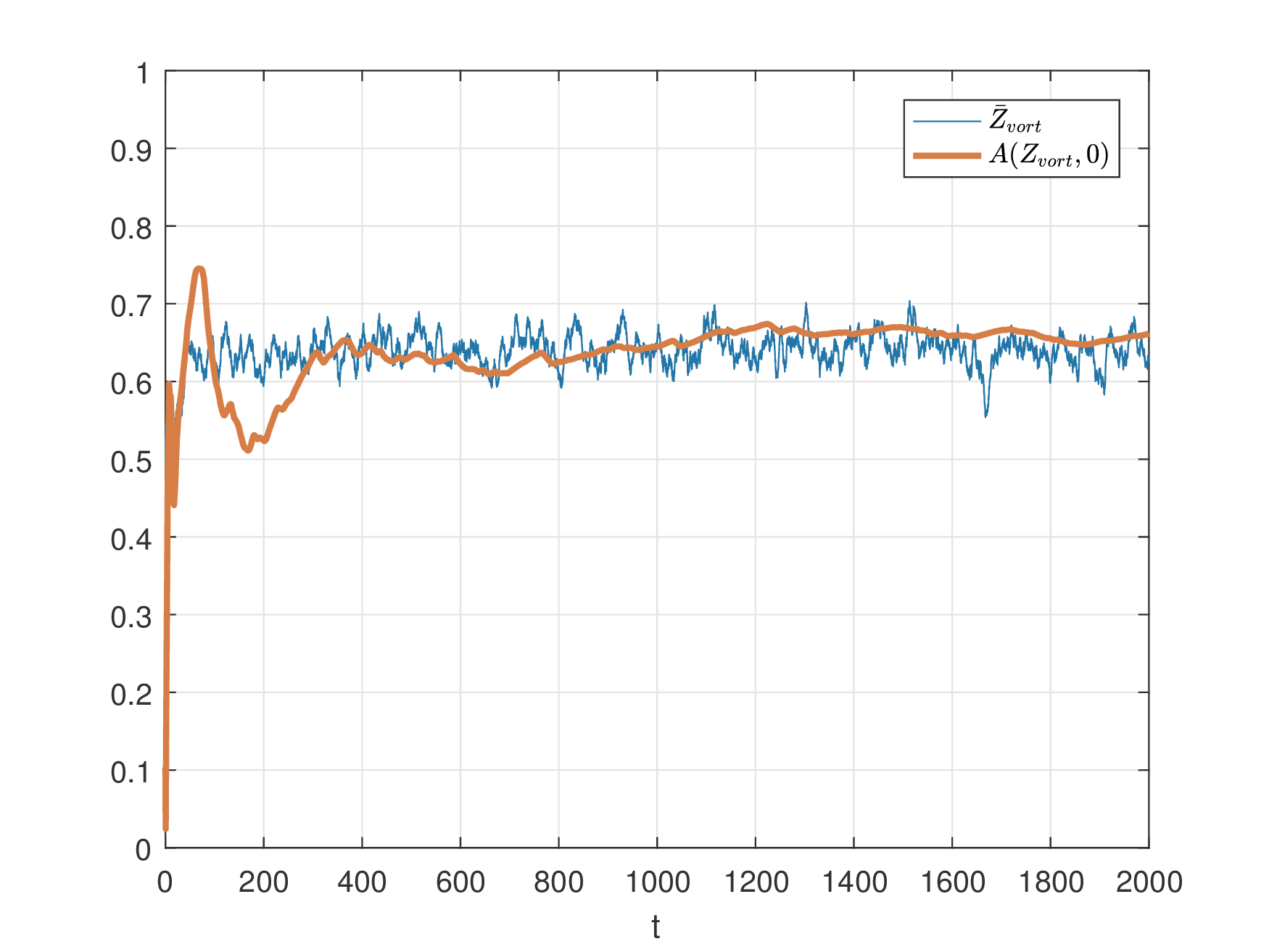}
  \caption{On average, the system is close to an $x$-bar state.}\label{fig: transition average Z}
\end{subfigure}
\caption{A single trajectory and the Monte Carlo average
  $\bar{Z}_{vort}(t)$ for $\delta=1.04$.}
\label{fig:delta1.04}
\end{figure}

For completeness, see also the discussion in \S\ref{S:Conclusion}, we
also include in Figure \ref{fig:delta1.04} a simulation with
$\nu=0.001$ that represents a single sample path for $\delta=1.04$, the value of $\delta$ for which transitions among the quasi-stationary states were observed in \cite{BouchetSimonnet09}.

\begin{figure}[H]
        \centering
        \begin{subfigure}[b]{0.32\textwidth}
            \centering
            \includegraphics[width=\textwidth]{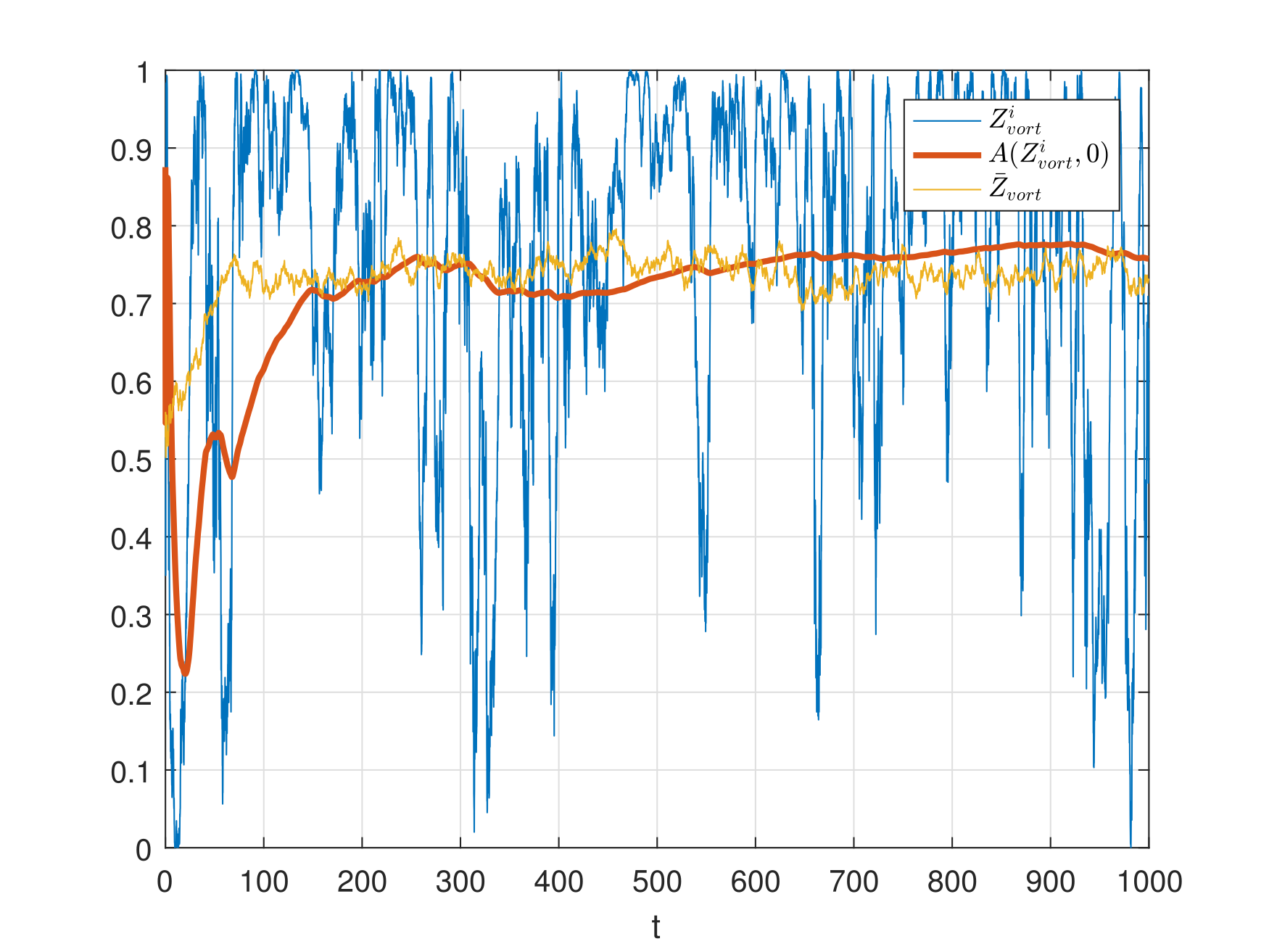}
            \caption{$\delta=1.10$}
            \label{fig:time avg 1.10}
        \end{subfigure}
        \begin{subfigure}[b]{0.32\textwidth}
            \includegraphics[width=\textwidth]{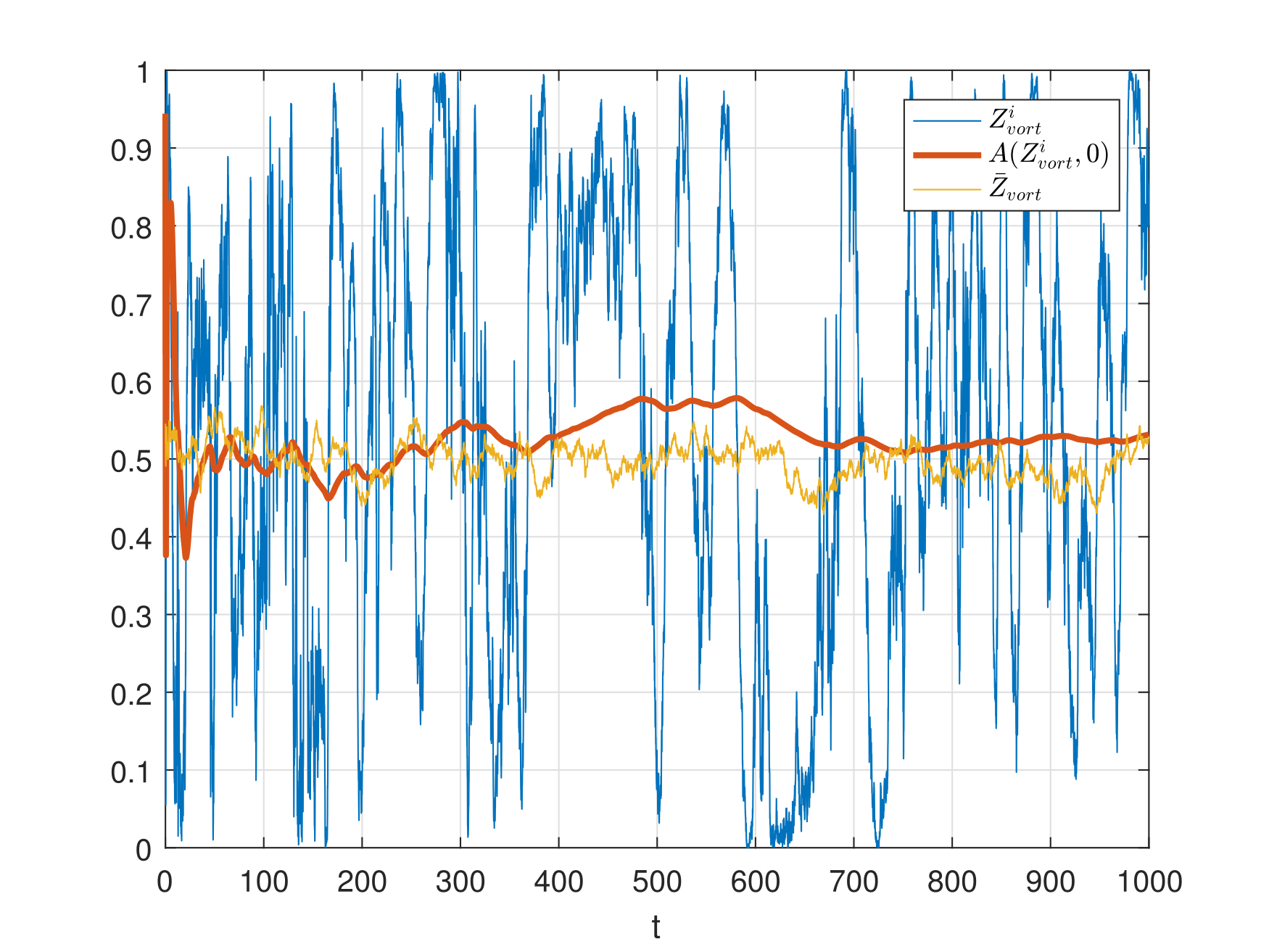}
            \caption{$\delta=1.0$}
            \label{fig:time avg 1.0}
        \end{subfigure}
          \begin{subfigure}[b]{0.32\textwidth}
            \includegraphics[width=\textwidth]{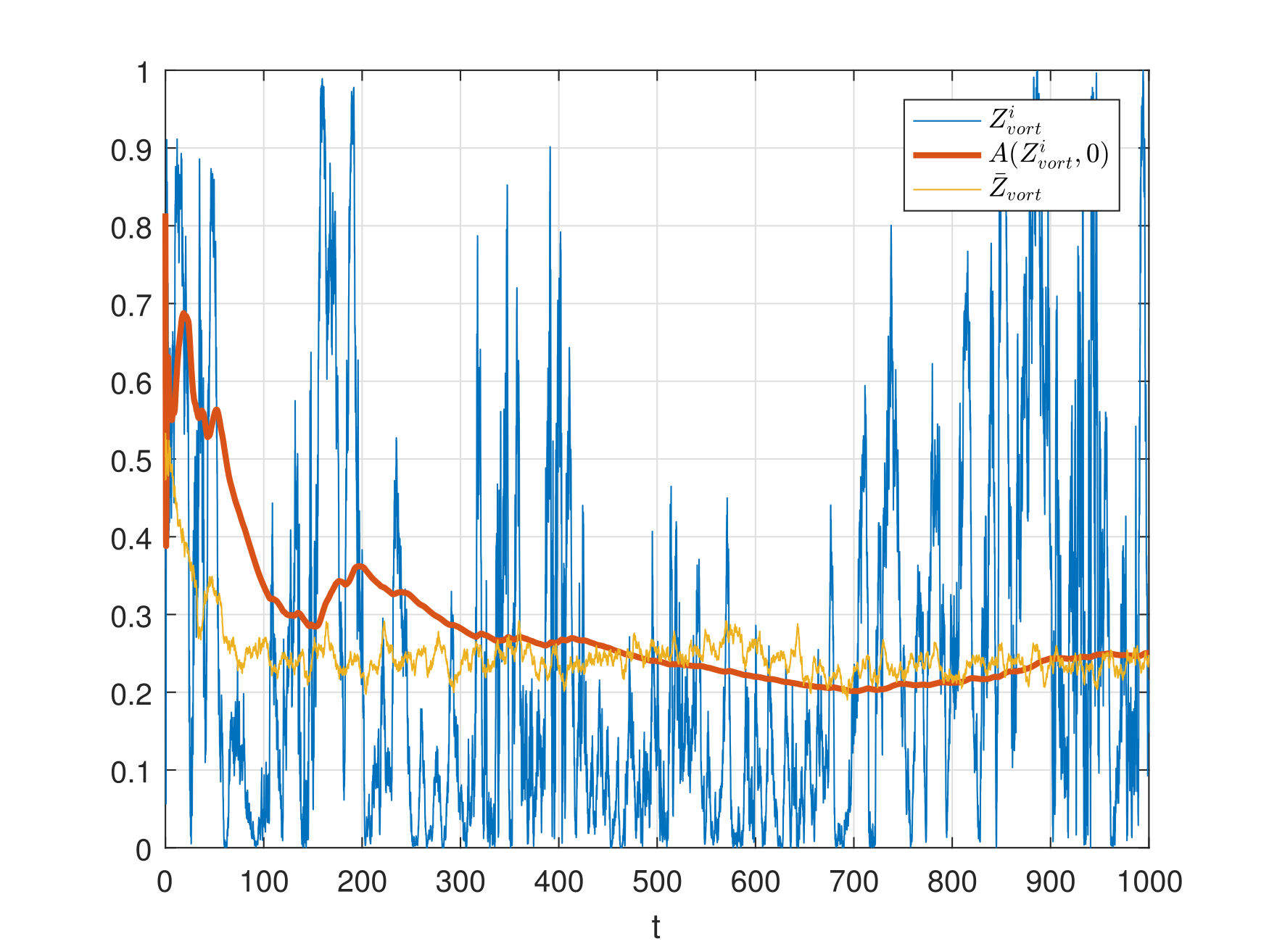}
            \caption{$\delta=0.90$}
            \label{fig:time avg 00.9}
        \end{subfigure}
                \caption{Comparing individual time average of a sample path with Monte Carlo average}
        \label{fig:time avg}
\end{figure}

Figure \ref{fig:sample_path} shows that individual trajectories exhibit transitions between quasi-stationary states, visiting the dipole and both bar states, as also observed in \cite{BouchetSimonnet09} for the same value of $\delta$. However Figure \ref{fig: transition average Z} shows that $\mathbb{E}[Z_{vort}(t)]$ picks the dominant state.
We now compute the time average of a randomly selected individual trial, given by $A(Z_{vort},t_{burn})$, to confirm that it tracks the Monte Carlo average, $\bar{Z}_{vort}(t)$. 
Figure \ref{fig:time avg} shows two things. First, for the given values of $\delta$,  a sample path may experience many transitions among the quasi-stationary states. Second, the time average of the sample path does eventually track the Monte Carlo average.

\subsection{Reduced Model}
We now turn our attention to the reduced model \eqref{E:deltanot1}. We confirm numerically that the reduced model captures the qualitative dynamics of the full vorticity equation with regard to the dominant quasi-stationary state.

As we mentioned in Remark \ref{R:RealSystem} we will be working with the real system in which $\omega_{1,3,5,7}$, as well as the Wiener processes, $W_{1,3,5,7}$, are all real valued. This leads to the following system which serves as the acting reduced model in the upcoming simulations.
\begin{equation}
    \begin{aligned}
        \dot \omega_1 &= - \frac{\nu}{\delta^2} \omega_1 + \frac{1}{\delta(1+\delta^2)}[\omega_3\omega_7 - \omega_3 \omega_5] + \frac{3\delta^6}{2\nu(4+\delta^2)(1+\delta^2)^2}\omega_1(\omega_5^2 + \omega_7^2) +\sqrt{2\nu}\sigma_{1}\dot{W}_{1} \\
        \dot \omega_3 &= - \nu \omega_3 + \frac{\delta^3}{(1+\delta^2)}[\omega_1\omega_5 - \omega_1 \omega_7] + \frac{3\delta^2}{2\nu(1+4\delta^2)(1+\delta^2)^2}\omega_3(\omega_5^2 + \omega_7^2) +\sqrt{2\nu}\sigma_{3}\dot{W}_{3}\label{E:real_system}  \\
        \dot \omega_5 &= - \nu \frac{1+\delta^2}{\delta^2} \omega_5 -\frac{\delta^2-1}{\delta}\omega_1 \omega_3  \\ & \quad  -\frac{\delta^6(3+\delta^2)}{2\nu(4+\delta^2)(1+\delta^2)} \omega_5 \omega_1^2-\frac{1+3\delta^2}{2\nu\delta^2(1+4\delta^2)(1+\delta^2)}\omega_5 \omega_3^2   +\sqrt{2\nu}\sigma_{5}\dot{W}_{5} \\
        \dot \omega_7 &= - \nu \frac{1+\delta^2}{\delta^2} \omega_7 +\frac{\delta^2-1}{\delta}\omega_1 \omega_3 \\ & \quad  -\frac{\delta^6(3+\delta^2)}{2\nu(4+\delta^2)(1+\delta^2)} \omega_7 \omega_1^2 - \frac{1+3\delta^2}{2\nu\delta^2(1+4\delta^2)(1+\delta^2)}\omega_7 \omega_3^2 +\sqrt{2\nu}\sigma_{7}\dot{W}_{7}.
    \end{aligned}
\end{equation}

To be consistent with the spatial decay of the noise in the simulations of the stochastically forced vorticity equation \eqref{E:vortS}, given by \eqref{E:noise-decay}, we choose
\[
\sigma_{1,3}=e^{-\alpha_0}\quad \mbox{and}\quad \sigma_{5,7}=e^{-2\alpha_0}.
\]

First we aim to establish that the reduced model \eqref{E:real_system} can serve as a good approximation to the vorticity equation with noise, \eqref{E:vortS}, for $\delta\approx1$. Second, it will be established that the selection of the bar or dipole state that dominates is consistent with the results of \cite{BeckCooperSpiliopoulos} for the deterministic equation: $x$-bar for $\delta>1$, $y$-bar for $\delta<1$, and dipole for $\delta=1$.

\begin{figure}[H]
	\centering
	\includegraphics[width=10cm, height=6cm]{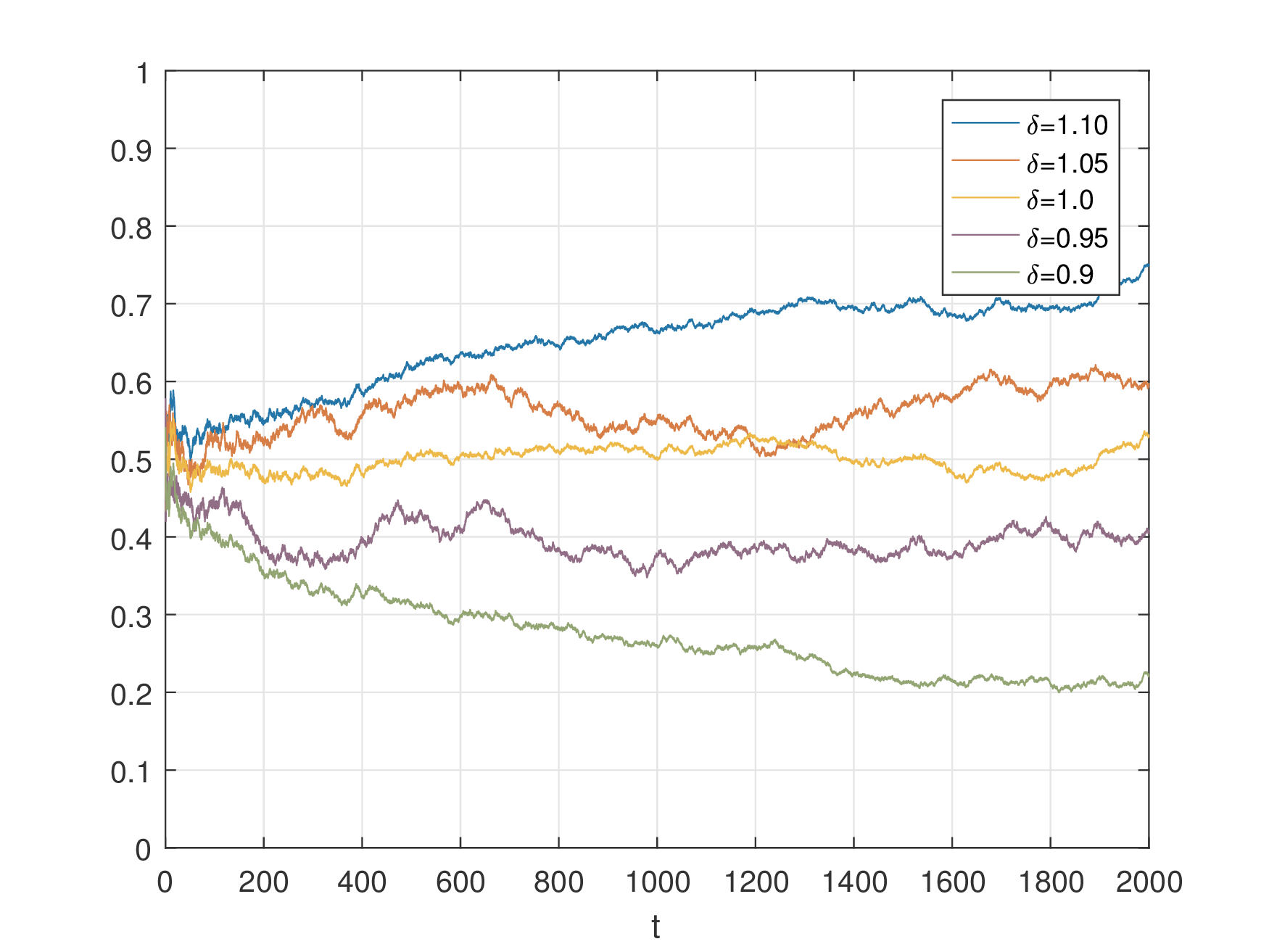}
	\caption{Simulation of $\bar{Z}_{red}(t)$ with noise for $\nu=0.001$.}\label{multidelta_scaled}
\end{figure}

Figure \ref{multidelta_scaled} shows numerical evidence supporting that the dynamics of the order parameter, governed by the reduced system \eqref{E:real_system}, follows the same trend as when the full vorticity equation is simulated.

\begin{figure}[H]
        \centering
        \begin{subfigure}[b]{0.32\textwidth}
            \centering
            \includegraphics[width=\textwidth]{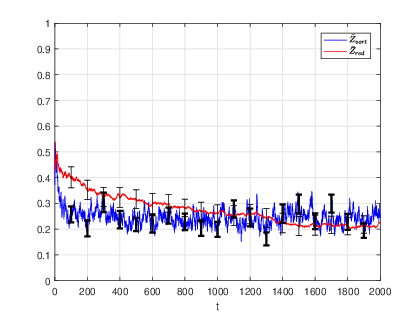}
            \caption{$\delta=0.90$}
            \label{fig:comp vort/red 90}
        \end{subfigure}
        \begin{subfigure}[b]{0.32\textwidth}
            \includegraphics[width=\textwidth]{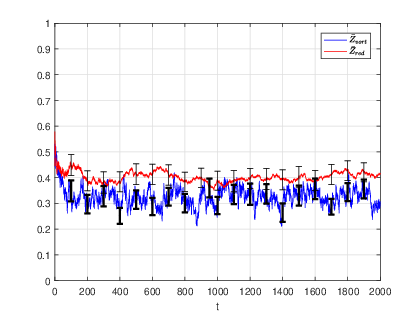}
            \caption{$\delta=0.95$}
            \label{fig:comp vort/red 95}
        \end{subfigure}
          \begin{subfigure}[b]{0.32\textwidth}
            \includegraphics[width=\textwidth]{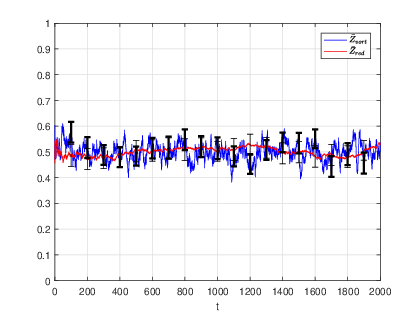}
            \caption{$\delta=1.0$}
            \label{fig:comp vort/red 1}
        \end{subfigure}

              \medskip
                \begin{subfigure}[b]{0.32\textwidth}
            \includegraphics[width=\textwidth]{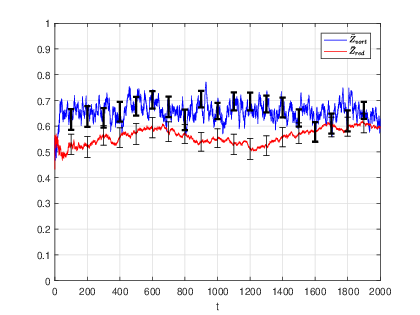}
            \caption{$\delta=1.05$}
            \label{fig:comp vort/red 105}
        \end{subfigure}
          \begin{subfigure}[b]{0.32\textwidth}
            \includegraphics[width=\textwidth]{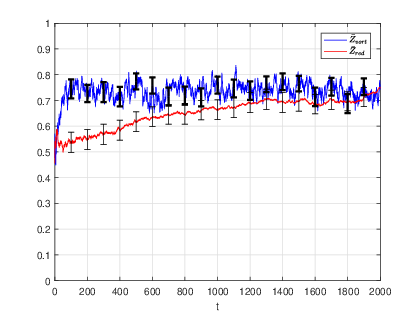}
            \caption{$\delta=1.10$}
            \label{fig:comp vort/red 110}
        \end{subfigure}
                \caption{Comparing $\bar{Z}_{vort}(t)$ and $\bar{Z}_{red}(t)$ averaged over N=200 trials with $\nu=0.001$. Corresponding 95 \% confidence error bars are also included.}

        \label{fig:comp vort/red}
\end{figure}

The plots of these Monte Carlo simulations (averaged over $N=200$ trials) show that the trend toward the appropriate quasi-stationary state is captured by the reduced model. Starting with zero initial conditions, when the noise is added, the simulations show that for $\delta>1$, the order parameter increases toward 1, indicating evolution to an $x$-bar state. Conversely, for $\delta<1$, the order parameter decreases toward a value corresponding to a $y$-bar state. Finally, when $\delta=1$, $\bar{Z}_{red}(t)$ remains near 1/2 indicating the system is in a dipole state. Figures \ref{fig:comp vort/red 90}-\ref{fig:comp vort/red 110} serve to compare the evolution of $\bar{Z}_{vort}(t)$ and $\bar{Z}_{red}(t)$
taken over $N=200$ trials, for values of $\delta$ close to 1. The bars denote the error for the 95\% confidence intervals for $\bar{Z}_{vort}$ (bold) and $\bar{Z}_{red}$ (thin).

One can see that $\bar{Z}_{vort}(t)$ and $\bar{Z}_{red}(t)$ both trend in the same direction, with similar variances (typically between 0.06-0.08 for $0\leq t\leq 2000$). Furthermore, their respective confidence intervals begin to converge until they overlap. Indeed the model can be used to determine towards which quasi-stationary state the system evolves for a given value of $\delta.$

\section{Perturbation Analysis}\label{S:PDE}

Motivated by the numerics from \S\ref{S:sims}, this section investigates the expected behavior of $Z_{red}(t)$ as $\delta\rightarrow 1$ while viewing the problem as a perturbation from the $\delta=1$  and $\nu=0$ case. Using the backward Kolmogorov equation associated to \eqref{E:real_system}, the goal is to derive a system of PDE that will provide insight on how the expected value of $Z_{red}(t)$, to leading order, depends on values of $\delta$ close to 1. To do this we pose the problem as a perturbation of the spatial domain, setting $\delta^2=1+\epsilon_0\epsilon$. Here, $0<\epsilon \ll 1$ acts as the small perturbation parameter and $\epsilon_0=\pm1$ determines which dimension of the torus is longer. Following known homogenization techniques, see for example \cite{PavliotisStewart},  we scale \eqref{E:real_system} in a way that reveals a slow-fast system of SDE. Then, we write the backward Kolmogorov equation to reach the ultimate goal of determining equations that govern the limiting evolution of $\mathbb{E}[Z_{red}(t)]$ as $\epsilon\rightarrow0$ once the fast variables are averaged out. First, for ease of notation, rename the dependent variables as follows,
\[
\tilde{p}:=\mbox{Re}(\omega_1), \quad \tilde{q}:=\mbox{Re}(\omega_3), \quad \tilde{r}:=\mbox{Re}(\omega_5), \quad \tilde{s}:=\mbox{Re}(\omega_7).
\]

Now \eqref{E:real_system} can be expressed as
\begin{eqnarray}
        \dot{\tilde{p}} &=& - \frac{\nu}{\delta^2} \tilde{p} + \frac{1}{\delta(1+\delta^2)}\tilde{q}(\tilde{s}-\tilde{r}) + \frac{3\delta^6}{2\nu(4+\delta^2)(1+\delta^2)^2}\tilde{p}(|\tilde{r}|^2 + |\tilde{s}|^2) +\sqrt{2\nu}\sigma_{1}\dot{W}_{1} \nonumber \\
        \dot{\tilde{q}} &=& - \nu \tilde{q} + \frac{\delta^3}{(1+\delta^2)}\tilde{p}(\tilde{r}-\tilde{s}) + \frac{3\delta^2}{2\nu(1+4\delta^2)(1+\delta^2)^2}\tilde{q}(|\tilde{r}|^2 + |\tilde{s}|^2) +\sqrt{2\nu}\sigma_{3}\dot{W}_{3}  \label{E:xy} \\
        \dot{\tilde{r}} &=& - \nu \frac{1+\delta^2}{\delta^2} \tilde{r} -\frac{\delta^2-1}{\delta}\tilde{p}\tilde{q}   -\frac{\delta^6(3+\delta^2)}{2\nu(4+\delta^2)(1+\delta^2)} \tilde{r} |\tilde{p}|^2-\frac{1+3\delta^2}{2\nu\delta^2(1+4\delta^2)(1+\delta^2)}\tilde{r} |\tilde{y}|^2   +\sqrt{2\nu}\sigma_{5}\dot{W}_{5} \nonumber \\
        \dot{\tilde{s}} &=& - \nu \frac{1+\delta^2}{\delta^2} \tilde{s} +\frac{\delta^2-1}{\delta}\tilde{p} \tilde{q} -\frac{\delta^6(3+\delta^2)}{2\nu(4+\delta^2)(1+\delta^2)} \tilde{s} |\tilde{p}|^2 - \frac{1+3\delta^2}{2\nu\delta^2(1+4\delta^2)(1+\delta^2)}\tilde{s} |\tilde{q}|^2 +\sqrt{2\nu}\sigma_{7}\dot{W}_{7}.\nonumber
\end{eqnarray}

Before inserting the Taylor expansions in $\epsilon$ for the
coefficients with $\delta^2=1+\epsilon_0\epsilon$, we first scale
\eqref{E:xy} appropriately to obtain a clear slow-fast system. As in
\cite{BeckCooperSpiliopoulos},  the low modes represented by
$\tilde{p}$ and $\tilde{q}$ correspond to the slow variables while the
high modes, $\tilde{r}$ and $\tilde{s}$, represent the fast
variables. Below, we give a more general version of the scaled
equations for just the $\tilde{p}$ (analogous to $\tilde{q}$) and
$\tilde{r}$ (analogous to $\tilde{s}$) equations. We use the following
space-time and parameter scalings: $\nu=\epsilon^\mu\nu_0$,
$\tilde{p}=\epsilon^\xi p$, $\tilde{q}=\epsilon^\xi q$,
$\tilde{r}=\epsilon^\eta r$, $\tilde{s}=\epsilon^\eta s$, and
$\tau=\epsilon^\gamma t$. To simplify the scaled equations, we will
relate $\mu$, $\xi$, $\eta$ and $\gamma$ to put the resulting system in a more desirable form. We neglect the $\epsilon$ dependence of $p$, $q$, $r$ and $s$ for readability. Below, the ``prime" notation denotes differentiation with respect to the scaled time variable, $\tau$.
\begin{equation*}
    \begin{aligned}
    p' &= \epsilon^{\mu-\gamma}(-\frac{\nu_0}{\delta^2}p)+\epsilon^{\eta-\gamma}\frac{1}{\delta(1+\delta^2)}q(s-r)+\epsilon^{2\eta-\mu-\gamma}\frac{3\delta^6}{2\nu_0(4+\delta^2)(1+\delta^2)^2}p(r^2+s^2) + \epsilon^{\frac{\mu-\gamma}{2}-\xi}\sqrt{2\nu_0}\sigma_1W'_5(\tau)  \\
    r' &= \epsilon^{\mu-\gamma}(-\nu_0\frac{1+\delta^2}{\delta^2}r)-\epsilon^{2\xi-\gamma-\eta}\frac{\delta^2-1}{\delta}(pq)-\epsilon^{2\xi-\mu-\gamma}\frac{1}{\nu_0}\left(\frac{\delta^6(3+\delta^2)}{2(4+\delta^2)(1+\delta^2)}rp^2+\frac{1+3\delta^2}{2(1+4\delta^2)(1+\delta^2)}rq^2\right) \\ & +\epsilon^{\frac{\mu-\gamma}{2}-\eta}\sqrt{2\nu_0}\sigma_5 W'(\tau)
    \end{aligned}
\end{equation*}

Now set $2\eta=\mu+\gamma$, $2\xi=\mu-\gamma \Rightarrow \gamma=\mu-2\xi$, with $0<\gamma<\xi<\frac{\mu}{2}<\eta<\mu.$ Then the fully scaled system (still neglecting Taylor expansions of $\delta$  in $\epsilon$ for now) becomes
\begin{eqnarray}
        p' &=& \epsilon^{2\xi}(- \frac{\nu_0}{\delta^2} p) + \epsilon^{\xi}\frac{1}{\delta(1+\delta^2)}q(s-r) + \frac{3\delta^6}{2\nu_0(4+\delta^2)(1+\delta^2)^2}p(|s|^2 + |s|^2) +\sqrt{2\nu_0}\sigma_{1}W_{1}'(\tau) \nonumber \\
        q' &=& \epsilon^{2\xi}(- \nu_0 q) + \epsilon^{\xi}\frac{\delta^3}{(1+\delta^2)}p(r-s) + \frac{3\delta^2}{2\nu_0(1+4\delta^2)(1+\delta^2)^2}q(|r|^2 + |s|^2) +\sqrt{2\nu_0}\sigma_{3}W'_{3}(\tau) \label{E:scaledSDE}  \\
        r' &=& \epsilon^{2\xi}\left(- \nu_0 \frac{1+\delta^2}{\delta^2}r\right) -\epsilon^{3\xi-2\eta}\frac{\delta^2-1}{\delta}pq   -\epsilon^{2(\xi-\eta)}\left(\frac{\delta^6(3+\delta^2)}{2\nu_0(4+\delta^2)(1+\delta^2)} r |p|^2+\frac{1+3\delta^2}{2\nu_0\delta^2(1+4\delta^2)(1+\delta^2)}r|q|^2\right) \nonumber \\   &&+\epsilon^{\xi-\eta}\sqrt{2\nu_0}\sigma_{5}W'_{5}(\tau) \nonumber  \\
        s' &=& \epsilon^{2\xi}\left(- \nu_0 \frac{1+\delta^2}{\delta^2}s\right) +\epsilon^{3\xi-2\eta}\frac{\delta^2-1}{\delta}pq -\epsilon^{2(\xi-\eta)}\left(\frac{\delta^6(3+\delta^2)}{2\nu_0(4+\delta^2)(1+\delta^2)} s |p|^2 - \frac{1+3\delta^2}{2\nu_0\delta^2(1+4\delta^2)(1+\delta^2)}s |q|^2\right) \nonumber \\ &&+\epsilon^{\xi-\eta}\sqrt{2\nu_0}\sigma_{7}W'_{7}(\tau). \nonumber
\end{eqnarray}

For the scaled SDE \eqref{E:scaledSDE}, let $b^\epsilon=(b_{p}^\epsilon,b_{q}^\epsilon,b_{r}^\epsilon,b_{s}^\epsilon)$ denote the drift vector and $\Sigma^\epsilon(p,q,r,s;\sigma_{1},\sigma_{3},\sigma_{5},\sigma_{7})$ denote the diffusion matrix so that \eqref{E:scaledSDE} can be expressed, for $\mathcal{X}^\epsilon=(p,q,r,s)$ and $\frac{d\mathcal{W}}{d\tau}=\left(\frac{dW_1}{d\tau}, \frac{dW_3}{d\tau}, \frac{dW_5}{d\tau},\frac{dW_7}{d\tau}\right)$, as
\begin{equation}\label{E:SDE}
\frac{d\mathcal{X^\epsilon}}{d\tau}=b^\epsilon+\Sigma^\epsilon \frac{d\mathcal{W}}{d\tau}.
\end{equation}

Now replacing the $\delta$ coefficients appearing in \eqref{E:scaledSDE} with their Taylor expansions for $\delta^2=1+\epsilon_0\epsilon$ up to $\mathcal{O}(\epsilon^3)$, the drift vector is given by (still supressing $\epsilon$ dependence of $p,q,r,s$),
\begin{eqnarray}
b_{p}^\epsilon &=& \frac{1}{\nu_0}\left(\frac{3}{40}+\epsilon\epsilon_0\frac{27}{200}+\epsilon^2\frac{117}{4000}-\epsilon_0\epsilon^3\frac{123}{5000}\right)p(r^2+s^2) \nonumber \\
& & +\epsilon^{\xi}\left(\frac{1}{2}-\frac{\epsilon_0}{2}\epsilon+\frac{7}{16}\epsilon^2-\frac{3\epsilon_0}{8}\epsilon^3\right)q(s-r)-\epsilon^{2\xi}\nu_0(1-\epsilon_0\epsilon+\epsilon^2-\epsilon_0\epsilon^3)p \label{E:drift} \\
b_q^\epsilon &=& \frac{1}{\nu_0}\left(\frac{3}{40}-\frac{3}{50}\epsilon_0\epsilon+\frac{117}{4000}\epsilon^2-\frac{93}{20000}\epsilon_0\epsilon^3\right)q(r^2+s^2) \nonumber \\
& & +\epsilon^\xi \left(\frac{1}{2}+\frac{\epsilon_0}{2}-\frac{1}{8}\epsilon^2\right)p(r-s) -\epsilon^{2\xi}\nu_0q \nonumber \\
b_r^\epsilon &=& \epsilon^{2(\xi-\eta)}\frac{1}{\nu_0}\left[-\frac{1}{5}r(p^2+q^2)-\epsilon_0\epsilon\frac{1}{100}r(51p^2-31q^2)-\epsilon^{2}\frac{373}{1000}r(p^2+q^2)-\epsilon_0\epsilon^{3}\frac{1}{10000}r\left(379p^2-4109q^2\right)\right] \nonumber \\
& & -\epsilon^{2\xi}\nu_0(2-\epsilon_0\epsilon+\epsilon^2-\epsilon_0\epsilon^3)r -\epsilon^{3\xi-2\eta}\left(\epsilon_0\epsilon-\frac{1}{2}\epsilon^2+\frac{3}{8}\epsilon_0\epsilon^3\right)pq \nonumber \\
b_s^\epsilon &=& \epsilon^{2(\xi-\eta)}\frac{1}{\nu_0}\left[-\frac{1}{5}s(p^2+q^2)-\epsilon_0\epsilon\frac{1}{100}s(51p^2-31q^2)-\epsilon^{2}\frac{373}{1000}s(p^2+q^2)-\epsilon_0\epsilon^{3}\frac{1}{10000}s\left(379p^2-4109q^2\right)\right] \nonumber \\
& &  -\epsilon^{2\xi}\nu_0(2-\epsilon_0\epsilon+\epsilon^2-\epsilon_0\epsilon^3)s +\epsilon^{3\xi-2\eta}\left(\epsilon_0\epsilon-\frac{1}{2}\epsilon^2+\frac{3}{8}\epsilon_0\epsilon^3\right)pq \nonumber
\end{eqnarray}
and the diffusion matrix by
\begin{align}
\Sigma^\epsilon(p,q,r,s;\sigma_{1},\sigma_{3},\sigma_{5},\sigma_{7})
&=
\left( {\begin{array}{cccc}
\sqrt{2\nu_0}\sigma_1  & 0 & 0  & 0 \\
 0 & \sqrt{2\nu_0}\sigma_3 &0  & 0 \\
0 & 0 & \epsilon^{\xi-\eta}\sqrt{2\nu_0}\sigma_5  & 0\\
0 & 0 & 0 & \epsilon^{\xi-\eta}\sqrt{2\nu_0}\sigma_7  \\
 \end{array} } \right).\label{E:diffusion}
\end{align}

With $\mbox{H}(u^{\epsilon})$ denoting the Hessian matrix of $u^{\epsilon}$, we now write the backward Kolmogorov equation for \eqref{E:SDE}, which is defined as
\begin{equation}
\begin{aligned}
\frac{\partial  u^{\epsilon}}{\partial \tau} = b^\epsilon \cdot \nabla u^\epsilon + \frac{1}{2}\mbox{Tr}[(\Sigma^{\epsilon})^2 \mbox{H}(u^\epsilon)], \qquad &\mbox{in } \mathbb{R}^4 \times [0,T] \label{E:kol} \\
u^\epsilon(p,q,r,s,0) = \phi(p,q), \qquad\qquad\qquad &\mbox{on } \mathbb{R}^4\times \{0\}.
\end{aligned}
\end{equation}

The backward Kolmogorov equation has the useful property that the evolution of $u^\epsilon(\mathcal{X}^\epsilon,\tau)$ gives
\[
u^\epsilon(p,q,r,s,\tau)=\mathbb{E}\left[\phi(p_{\tau},q_{\tau}) \middle| \mbox{ } p_\tau(0)=p, \mbox{ } q_\tau(0)=q, \mbox{ } r_\tau(0)=r, \mbox{ } s_\tau(0)=s \right].
\]

Thus one ultimately is interested in initializing \eqref{E:kol} with $\phi(p,q)=Z_{red}=\frac{p^2}{p^2+q^2}$, but for now we proceed with a general initial condition, $\phi$. We seek a solution to \eqref{E:kol} that takes the form
\begin{equation}\label{E:u_eps}
    u^{\epsilon}(p,q,r,s,\tau)=u_0(p,q,r,s,\tau)+\epsilon u_1(p,q,r,s,\tau)+\epsilon^2 u_2(p,q,r,s,\tau)+\dots
\end{equation}
and wish to find the limiting dynamics, $u^\epsilon$ as $\epsilon\to 0$. As the goal is to identify the leading order expansion for $u^\epsilon$ we  determine a system of PDEs for $u_0$, $u_1$, and $u_2$.  We present now the calculations that lead to the characterization of  $u_0$, $u_1$, and $u_2$, see (\ref{E:u0u1new}).

Define $\mathcal{L}^\epsilon$ to be the operator acting on the right hand side of \eqref{E:kol}, so that $\frac{\partial u^\epsilon}{\partial \tau}=\mathcal{L}^\epsilon u^\epsilon$. Decomposing $\mathcal{L}^\epsilon$ by powers of $\epsilon$ using the expressions for $b^\epsilon$ and $\Sigma^\epsilon$ given in \eqref{E:drift} and \eqref{E:diffusion}, we write
\begin{eqnarray*}\label{E:kol_scaled}
\mathcal{L}^\epsilon u^\epsilon &=& \epsilon^{2(\xi-\eta)}\mathcal{L}_0u^\epsilon + \epsilon^{2(\xi-\eta)+1}\mathcal{L}_1u^\epsilon + \epsilon^{2(\xi-\eta)+2}\mathcal{L}_2 u^\epsilon+\epsilon^{2(\xi-\eta)+3}\mathcal{L}_3 u^\epsilon \nonumber \\
&+& \mathcal{L}_{4}u^\epsilon + \epsilon\mathcal{L}_{5} u^\epsilon+ \epsilon^2\mathcal{L}_6+\epsilon^3\mathcal{L}_7u^\epsilon \nonumber \\
&+&\epsilon^{\xi}\mathcal{L}_{8}u^\epsilon + \epsilon^{\xi+1}\mathcal{L}_{9}u^\epsilon + \epsilon^{\xi+2}\mathcal{L}_{10}u^\epsilon + \epsilon^{\xi+3}\mathcal{L}_{11}u^\epsilon \nonumber  \\
&+& \epsilon^{2\xi}\mathcal{L}_{12}u^\epsilon +\epsilon^{2\xi+1}\mathcal{L}_{13}u^\epsilon+ \epsilon^{2\xi+2}\mathcal{L}_{14}u^\epsilon+ \epsilon^{2\xi+3}\mathcal{L}_{15}u^\epsilon \nonumber  \\
&+& \epsilon^{3\xi-2\eta+1}\mathcal{L}_{16}u^\epsilon +\epsilon^{3\xi-2\eta+2}\mathcal{L}_{17}u^\epsilon +\epsilon^{3\xi-2\eta+3}\mathcal{L}_{18}u^\epsilon.
\end{eqnarray*}

A select few of the operators, $\mathcal{L}_i$, $i=1,\dots,19$,  that are most important in computing the leading order equations is provided in \eqref{E:operators}. A complete list of the expressions of the 19 operators can be found in the appendix.
\begin{eqnarray}
\mathcal{L}_0u &=& -\frac{1}{5\nu_0}(p^2+q^2)\left(r\frac{\partial u}{\partial r} + s\frac{\partial u}{\partial s}\right)+2\nu_0\left(\sigma_5^2\frac{\partial^2 u}{\partial r^2}+\sigma_7^2  \frac{\partial^2 u}{\partial s^2}\right)\nonumber \\
\mathcal{L}_1u &=& -\frac{\epsilon_0}{100\nu_0}(51p^2-31q^2)\left(r\frac{\partial u}{\partial r} + s\frac{\partial u}{\partial s}\right) \nonumber \\
\mathcal{L}_{4}u &=& \frac{3}{40\nu_0}(r^2+s^2)\left(p\frac{\partial u}{\partial p}+q\frac{\partial u}{\partial q}\right)+2\nu_0\left(\sigma_1^2\frac{\partial^2}{\partial p^2}+\sigma_3^2\frac{\partial^2}{\partial q^2}\right) \label{E:operators} \\
\mathcal{L}_{5}u  &=&\frac{\epsilon_0}{\nu_0}(r^2+s^2)\left(\frac{27}{200}p\frac{\partial u}{\partial p}-\frac{3}{50}q\frac{\partial u}{\partial q}\right) \nonumber \\\nonumber \\
\mathcal{L}_{8}u &=& -\frac{1}{2}(r-s)\left(q\frac{\partial u}{\partial p}-p\frac{\partial u}{\partial q}\right) \nonumber
\end{eqnarray}

Now we choose explicit values for $\xi$ and $\eta$ to obtain a simpler, but still representative system: $\eta=2$, $\xi=1$, hence $\mu=3$ ($\nu=\epsilon^\mu \nu_0=\epsilon^3\nu_0$). Note that this choice of parameters is not unique. Our goal is to make a choice that makes the computations that follow tractable, preserves a clear slow-fast system in the scaled equations \eqref{E:scaledSDE} and results in a system that exhibits qualitatively the same behavior as the original system in terms of the selection mechanism to $x$-bar, $y$-bar and dipole states depending on the values of $\delta$.   With this choice, the generator $\mathcal{L}^\epsilon$ of the PDE becomes
\begin{eqnarray}\label{E:simple_kol}
\mathcal{L}^\epsilon u^{\epsilon}&=&\epsilon^{-2}\mathcal{L}_0u^{\epsilon}+\epsilon^{-1}\mathcal{L}_1u^{\epsilon}+\left(\mathcal{L}_2+\mathcal{L}_{4}+\mathcal{L}_{16}\right)u^{\epsilon}+\epsilon\left(\mathcal{L}_3+\mathcal{L}_{5}+\mathcal{L}_{8}+\mathcal{L}_{8}\right) u^{\epsilon} \\
&+&\epsilon^2\left(\mathcal{L}_6+\mathcal{L}_{9}+\mathcal{L}_{12}+\mathcal{L}_{18}\right)u^{\epsilon}+\epsilon^3\left(\mathcal{L}_7+\mathcal{L}_{10}+\mathcal{L}_{13}\right)u^{\epsilon} +\epsilon^4\left(\mathcal{L}_{11}+\mathcal{L}_{14}\right)u^{\epsilon}+\epsilon^5\mathcal{L}_{15}u^{\epsilon}\nonumber.
\end{eqnarray}

The ansatz given in \eqref{E:u_eps} can now be inserted into the backward Kolmogorov equation \eqref{E:kol} using the expression of $\mathcal{L}^\epsilon$ given above in \eqref{E:simple_kol}. Matching coefficients on both sides of the equation yields the following leading order equations,
\begin{subequations}
\begin{eqnarray}
&\mathcal{O}(\epsilon^{-2}):& -\mathcal{L}_0u_0=0  \label{E:init_beta2a} \\
&\mathcal{O}(\epsilon^{-1}):& -\mathcal{L}_0u_1=\mathcal{L}_1u_0 \label{E:init_beta2b}  \\
&\mathcal{O}(1):& -\mathcal{L}_0u_2=-\frac{\partial u_0}{\partial\tau}+\mathcal{L}_1u_1+\left(\mathcal{L}_2+\mathcal{L}_{4}+\mathcal{L}_{16}\right)u_0\label{E:init_beta2c} \\
&\mathcal{O}(\epsilon):& -\mathcal{L}_0u_3=-\frac{\partial u_1}{\partial\tau}+\mathcal{L}_1u_2+(\mathcal{L}_2+\mathcal{L}_{4}+\mathcal{L}_{16})u_1+\left(\mathcal{L}_3+\mathcal{L}_{5}+\mathcal{L}_{8}+\mathcal{L}_{17}\right)u_0\label{E:init_beta2d}.
\end{eqnarray}
\end{subequations}

Equation (\ref{E:init_beta2a}) implies $u_0$ lies in the kernel of $\mathcal{L}_0$, which elliptic PDE theory tells us contains only functions constant in $r$ and $s$. Since $\mathcal{L}_1$ is also a differential operator in $r$ and $s$ only, \eqref{E:init_beta2b} implies that $u_1$ is constant in $r$ and $s$ as well. One can also see that $u_0$ and $u_1$ are in the kernel of each of $\mathcal{L}_{2,3,16,17}$ (see appendix). Hence the leading order system given by \eqref{E:init_beta2a}-\eqref{E:init_beta2d} can be reduced to
\begin{subequations}
\begin{eqnarray}
&\mathcal{O}(\epsilon^{-2}):& -\mathcal{L}_0u_0=0 \Rightarrow u_0=u_0(p,q,\tau) \label{E:beta2a} \\
&\mathcal{O}(\epsilon^{-1}):& -\mathcal{L}_0u_1=\mathcal{L}_1u_0 \Rightarrow u_1=u_1(p,q,\tau)\label{E:beta2b}  \\
&\mathcal{O}(1):& -\mathcal{L}_0u_2=-\frac{\partial u_0}{\partial\tau}+\mathcal{L}_{4}u_0\label{E:beta2c} \\
&\mathcal{O}(\epsilon):& -\mathcal{L}_0u_3=-\frac{\partial u_1}{\partial\tau}+\mathcal{L}_1u_2+\mathcal{L}_{4}u_1+\left(\mathcal{L}_{5}+\mathcal{L}_{8}\right)u_0,\label{E:beta2d}
\end{eqnarray}
\end{subequations}
where $\mathcal{L}_0$, $\mathcal{L}_1$, $\mathcal{L}_4$, $\mathcal{L}_5$, and $\mathcal{L}_8$ are presented in \eqref{E:operators}. Let $\rho^\infty(r,s;p,q)$ be the stationary density that satisfies the adjoint problem
\[
\mathcal{L}^{*}_{0} \rho^{\infty}(r,s;p,q)=0.
\]

Once $\rho^\infty$ is known we can integrate against the invariant measure to obtain the solvability conditions for equations (\ref{E:beta2c}) and (\ref{E:beta2d})
\begin{eqnarray}
\frac{\partial u_0}{\partial\tau}&=&\int_{\mathbb{R}^2} \mathcal{L}_{4}u_0 \rho^{\infty}(r,s;p,q)\rmd r \rmd s \label{E:solvability} \\
\frac{\partial u_1}{\partial\tau}&=&\int_{\mathbb{R}^2} \left(\mathcal{L}_1u_2+\mathcal{L}_{4}u_1+\left(\mathcal{L}_{5}+\mathcal{L}_{8}\right)u_0\right) \rho^{\infty}(r,s;p,q)\rmd r \rmd s. \nonumber
\end{eqnarray}

Before we consider the integrals in \eqref{E:solvability}, $\rho^{\infty}$ must be identified. The operator,
\[
\mathcal{L}_0=\left(-\frac{1}{5\nu_0}r(p^2+q^2),-\frac{1}{5\nu_0}s(p^2+q^2)\right)\cdot\left(\frac{\partial}{\partial r},\frac{\partial}{\partial s}\right)+\frac{1}{2}  \left( {\begin{array}{cc}
   4\nu_0\sigma_5^2   & 0  \\
   0 & 4\nu_0\sigma_7^2
  \end{array} } \right)\left({\begin{array}{c} \frac{\partial^2}{\partial r^2} \\ \frac{\partial^2}{\partial s^2}\end{array}} \right),
\]
corresponds to the backward Kolmogorov equation for the following system, parameterized by the fixed (slow) variables $p$ and $q$.
\begin{eqnarray}
\dot{\hat{r}}&=& -\frac{1}{5\nu_0}(p^2+q^2)\hat{r}+\sigma_5\sqrt{2\nu_0}\dot{W}_5 \nonumber \\
\dot{\hat{s}}&=& -\frac{1}{5\nu_0}(p^2+q^2)\hat{s}+\sigma_7 \sqrt{2\nu_0}\dot{W}_7. \nonumber
\end{eqnarray}

These processes are independent Ornstein-Uhlenbeck processes and are therefore Gaussian. The equilibrium (stationary) density which corresponds to $\rho^{\infty}(r,s;p,q)$ is that of the bivariate Gaussian distribution with
\[
r\sim N\left(0,\frac{5\nu_0^2\sigma_5^2}{p^2+q^2}\right), \qquad s\sim N\left(0,\frac{5\nu_0^2\sigma_7^2}{p^2+q^2}\right).
\]

Therefore the invariant joint density is
\[
\rho^\infty(r,s,p,q)=\frac{p^2+q^2}{10\pi\nu_0^2\sigma_5\sigma_7}e^{-\frac{p^2+q^2}{10\nu_0^2}(\frac{r^2}{\sigma_5^2}+\frac{s^2}{\sigma_7^2})}.
\]

To aid in the computations of the integrals given in (\ref{E:solvability}), the following integral evaluations will be useful and can be simply obtained through the mean and variance of the stationary distribution.
\begin{subequations}
\begin{eqnarray}
\int_{\mathbb{R}^2} (r^2+s^2)\rho^\infty(r,s;p,q)\rmd r\rmd s &=& \frac{5\nu_0^2}{p^2+q^2}(\sigma_5^2+\sigma_7^2) \label{E:integralA} \\
\int_{\mathbb{R}^2} (r-s)\rho^\infty(r,s;p,q)\rmd r\rmd s  &=& 0. \label{E:integralB}
\end{eqnarray}
\end{subequations}

Next consider the solvability conditions (\ref{E:solvability}) one at a time. From the $u_0$ equation and the integral \eqref{E:integralA},
\begin{eqnarray}
\frac{\partial u_0}{\partial \tau} &=& \int_{\mathbb{R}^2} \mathcal{L}_{4}u_0 \rho^\infty \rmd r\rmd s \nonumber \\
&=& 2\nu_0\left(\sigma_1^2\frac{\partial^2 u_0}{\partial p^2}+\sigma_3^2\frac{\partial^2 u_0}{\partial q^2}\right)+\frac{3}{40\nu_0}\left(p\frac{\partial u_0}{\partial p}+q\frac{\partial u_0}{\partial q}\right)\int_{\mathbb{R}^2}(r^2+s^2)\rho^\infty(r,s;p,q)\rmd r\rmd s \nonumber \\
&=& 2\nu_0\left(\sigma_1^2\frac{\partial^2 u_0}{\partial p^2}+\sigma_3^2\frac{\partial^2 u_0}{\partial q^2}\right)+\frac{3\nu_0}{8}(\sigma_5^2+\sigma_7^2)\left(\frac{p}{p^2+q^2}\frac{\partial u_0}{\partial p}+\frac{q}{p^2+q^2}\frac{\partial u_0}{\partial q}\right). \nonumber
\end{eqnarray}

From this we obtain the effective equations for $p$ and $q$ for small $\epsilon$ after the fast variables $r$ and $s$ are averaged out. The slow motion can hence be approximated, for $0<\epsilon\ll1$, by $\bar{p}$ and $\bar{q}$ governed by,
\begin{eqnarray}
\bar{p}' &=& \frac{3\nu_0}{8}(\sigma_5^2+\sigma_7^2)\frac{\bar{p}}{\bar{p}^2+\bar{q}^2} +\sigma_1\sqrt{2\nu_0}W_1' \nonumber \\
\bar{q}' &=& \frac{3\nu_0}{8}(\sigma_5^2+\sigma_7^2)\frac{\bar{q}}{\bar{p}^2+\bar{q}^2}+\sigma_3\sqrt{2\nu_0}W_3'. \label{E:effectiveslow}
\end{eqnarray}

Since $\epsilon_0$ dependence does not appear in the first order equations, we will need to determine $u_1$ to see its effects. Consider the solvability condition for $u_1$ in \eqref{E:solvability}. Computing this integral requires us to evaluate the following integrals.
\begin{eqnarray}
I_{0}&=&\int_{\mathbb{R}^2}\mathcal{L}_{5}u_0 \rho^\infty \rmd r\rmd s \nonumber \\
I_{0'} &=& \int_{\mathbb{R}^2} \mathcal{L}_{8}u_0\rho^\infty \rmd r \rmd s \nonumber \\
I_1 &=&\int_{\mathbb{R}^2} \mathcal{L}_{4}u_1\rho^\infty \rmd r \rmd s \nonumber \\
I_2 &=& \int_{\mathbb{R}^2} \mathcal{L}_1u_2\rho^\infty \rmd r\rmd s. \nonumber
\end{eqnarray}

In evaluating these, we see
\begin{eqnarray}
I_{0}&=&\int_{\mathbb{R}^2}  \frac{\epsilon_0}{\nu_0}\left(\frac{27}{200}p\frac{\partial u_0}{\partial p}-q\frac{\partial u_0}{\partial q}\right)(r^2+s^2) \rho^{\infty}(r,s;p,q) \rmd r\rmd s = \nonumber \\
&=& 5\epsilon_0\nu_0(\sigma_5^2+\sigma_7^2)\left(\frac{27}{200}\frac{p}{p^2+q^2}\frac{\partial u_0}{\partial p}-\frac{q}{p^2+q^2}\frac{\partial u_0}{\partial q}\right) \nonumber \\
I_{0'} &=& \int_{\mathbb{R}^2} -\frac{1}{2}\left(q\frac{\partial u_0}{\partial p}-p\frac{\partial u_0}{\partial q}\right)(r-s)\rho^\infty(r,s;p,q)\rmd r\rmd s=0 \nonumber \\
I_1 &=& 2\nu_0\left(\sigma_1^2\frac{\partial^2 u_1}{\partial p^2}+\sigma_3^2\frac{\partial^2 u_1}{\partial q^2}\right) \nonumber \\
& & + \frac{3\nu_0}{8}(\sigma_5^2+\sigma_7^2)\left(\frac{p}{p^2+q^2}\frac{\partial u_1}{\partial p}+\frac{q}{p^2+q^2}\frac{\partial u_1}{\partial q}\right) \nonumber \\
I_2 &=& -\epsilon_0\frac{1}{100\nu_0}(51p^2-31q^2)\int_{\mathbb{R}^2}\left(r\frac{\partial u_2}{\partial r}+s\frac{\partial u_2}{\partial s}\right)\rho^\infty(r,s;p,q) \rmd r\rmd s.
\end{eqnarray}

Since $u_2$ depends on the fast variables, this final integral cannot yet be computed. It will eventually be handled numerically.  Thus, formally, we have $u^\epsilon=u_0+\epsilon u_1+\epsilon^2 u_2+\dots$ satisfying,
\begin{eqnarray}
\frac{\partial u_0}{\partial \tau}&=&\frac{3\nu_0}{8}(\sigma_5^2+\sigma_7^2)\left(\frac{p}{p^2+q^2}\frac{\partial u_0}{\partial p}+\frac{q}{p^2+q^2}\frac{\partial u_0}{\partial q}\right)+2\nu_0\left(\sigma_1^2\frac{\partial^2 u_0}{\partial p^2}+\sigma_3^2\frac{\partial^2 u_0}{\partial q^2}\right) \label{E:u0u1new} \\
\frac{\partial u_1}{\partial \tau} &=&  \frac{3\nu_0}{8}(\sigma_5^2+\sigma_7^2)\left(\frac{p}{p^2+q^2}\frac{\partial u_1}{\partial p}+\frac{q}{p^2+q^2}\frac{\partial u_1}{\partial q}\right)+2\nu_0\left(\sigma_1^2\frac{\partial^2 u_1}{\partial p^2}+\sigma_3^2\frac{\partial^2 u_1}{\partial q^2}\right) \nonumber \\
&+& 5\epsilon_0\nu_0(\sigma_5^2+\sigma_7^2)\left(\frac{27}{200}\frac{p}{p^2+q^2}\frac{\partial u_0}{\partial p}-\frac{q}{p^2+q^2}\frac{\partial u_0}{\partial q}\right)-\frac{\epsilon_0}{100\nu_0}(51p^2-31q^2)\int_{\mathbb{R}^2} \left(r\frac{\partial u_2}{\partial r}+s\frac{\partial u_2}{\partial s}\right)\rmd \rho^\infty \nonumber \\
-\mathcal{L}_0u_2 &=& \frac{\partial u_0}{\partial \tau}-\mathcal{L}_{4}u_0\nonumber
\end{eqnarray}

We shall consider the system (\ref{E:u0u1new}) together with the initial conditions
\[
u_{0}(p,q,0)=\phi(p,q), \quad u_{1}(p,q,0)=0, \textrm{ and } u_{2}(p,q,r,s,0)=0.
\]

The PDE for $u_0$ immediately stands out as the backward Kolmogorov equation corresponding to the system given in (\ref{E:effectiveslow}). Despite its simple looking form, the regularity at the origin of the coefficients on the first derivative terms turn out to be a borderline case with regards to well posedness, see for example Chapter III, Section 1 of \cite{Ladyzenskaja}.  Nevertheless, we proceed formally and solve for $u_0$, $u_1$, and $u_2$ numerically after providing the initial condition
\[
\phi(p,q,0)=\frac{p^2}{p^2+q^2}
\]
so that $u^\epsilon(p,q,r,s,\tau)=\mathbb{E}[Z_{red}(\tau) | p_0=p, \mbox{ } q_0=q, \mbox{ } r_0=r, \mbox{ } s_0=s]$. The simulations of the system \eqref{E:u0u1new} provide an approximation to the deterministic evolution of
$\mathbb{E}[Z_{red}(\tau)]$ after averaging out the fast motion. The simulations of the system \eqref{E:u0u1new} provided in this section were conducted via finite differences on the domain $(p,q,r,s,\tau)\in[-5,5]^4\times[0,T]$ with Neumann boundary conditions.

\begin{figure}[H]
	\centering
	\includegraphics[width=10cm, height=6cm]{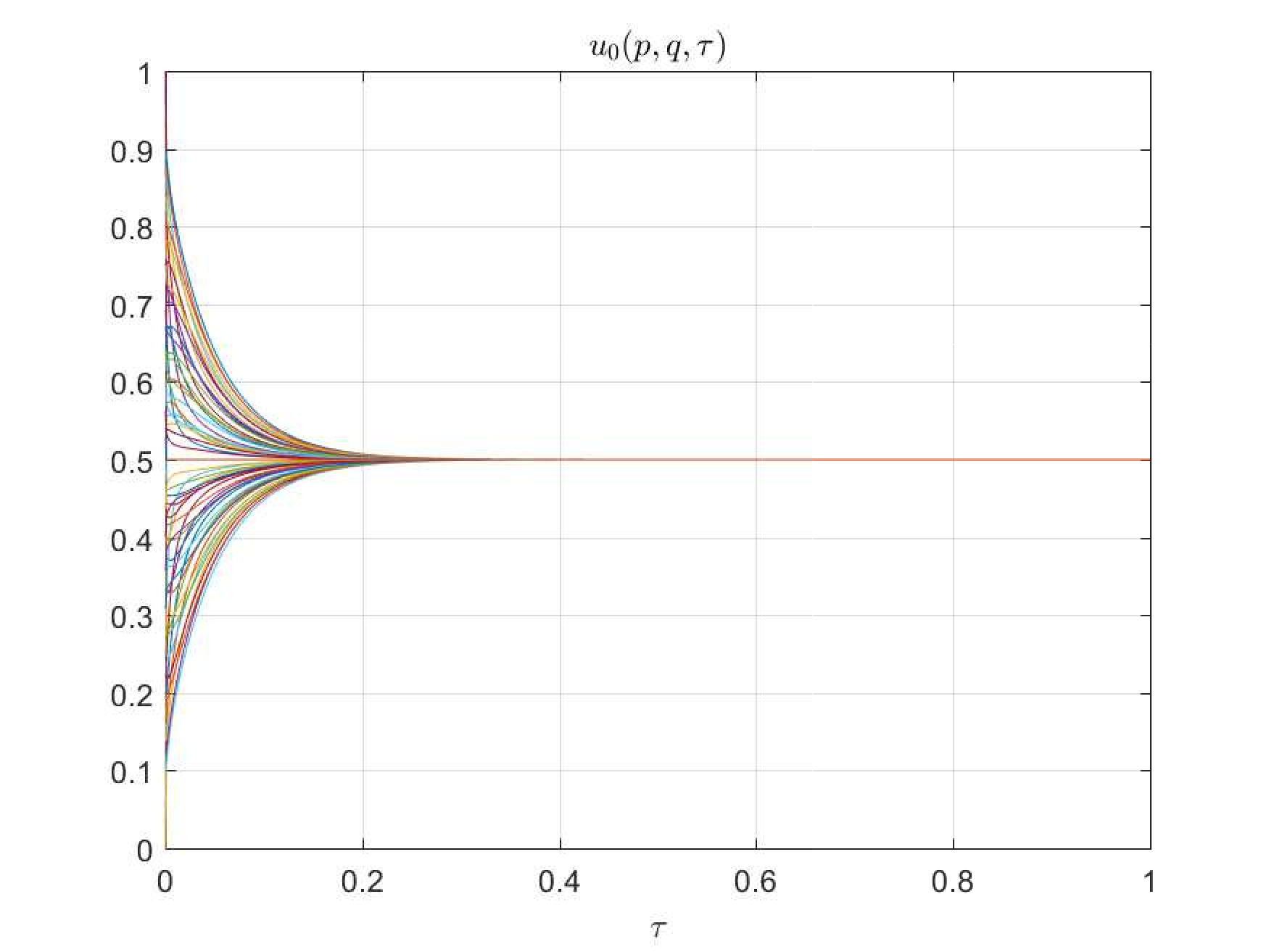}
	\caption{For $\epsilon=0$, $\mathbb{E}[Z_{red}(\tau)]\to1/2$}\label{fig: u0}
\end{figure}

When $\epsilon=0$, which implies $\delta=1$, the system is unperturbed and $u_0(p,q,t)=\mathbb{E}[Z_{red}(\tau) | p_0=p, \mbox{ } q_0=q]$. In this case, the results of the simulation show that the expected value of the order parameter $Z_{red}(\tau)$ converges to 1/2 for any initial values $p_0=p$ and $q_0=q$ of \eqref{E:deltanot1}, independent of $r$ and $s$. This indicates that the unperturbed system evolves to a dipole state, even if the initial state is close to an $x$- or $y$-bar state. Figure \ref{fig: u0} illustrates the evolution to a dipole for $\epsilon=0$ for several initial conditions $(p,q)$ chosen within the domain.

We conclude this section by comparing the numerical approximation to $\mathbb{E}[Z_{red}(\tau)]$, given by $u^\epsilon(\tau)$ to the average path of the order parameter $Z_{red}(t)$, i.e. $\bar{Z}_{red}(t)$, as computed via Monte Carlo simulation in \S\ref{S:sims}. Since the two models evolve on different timescales, we rescale $\tau$ so that our averaged PDE model is evolving on the original timescale. As such, suppressing the spatial $(p,q,r,s)$ dependency, let $\hat{u}^\epsilon(t):=u_0(t/\epsilon)+\epsilon u_1(t/\epsilon)=u_0(\tau)+\epsilon u_1(\tau)$ denote the $\mathcal{O}(\epsilon)$ approximation to $u^\epsilon$ in the original timescale, and let $\bar{Z}_{red}(t)$ denote the Monte Carlo average path of the order parameter under the dynamics of \eqref{E:deltanot1} obtained via the Monte Carlo simulations described in \S\ref{S:sims}. 

\begin{figure}[H]
        \centering
        \begin{subfigure}[b]{0.475\textwidth}
            \centering
            \includegraphics[width=\textwidth]{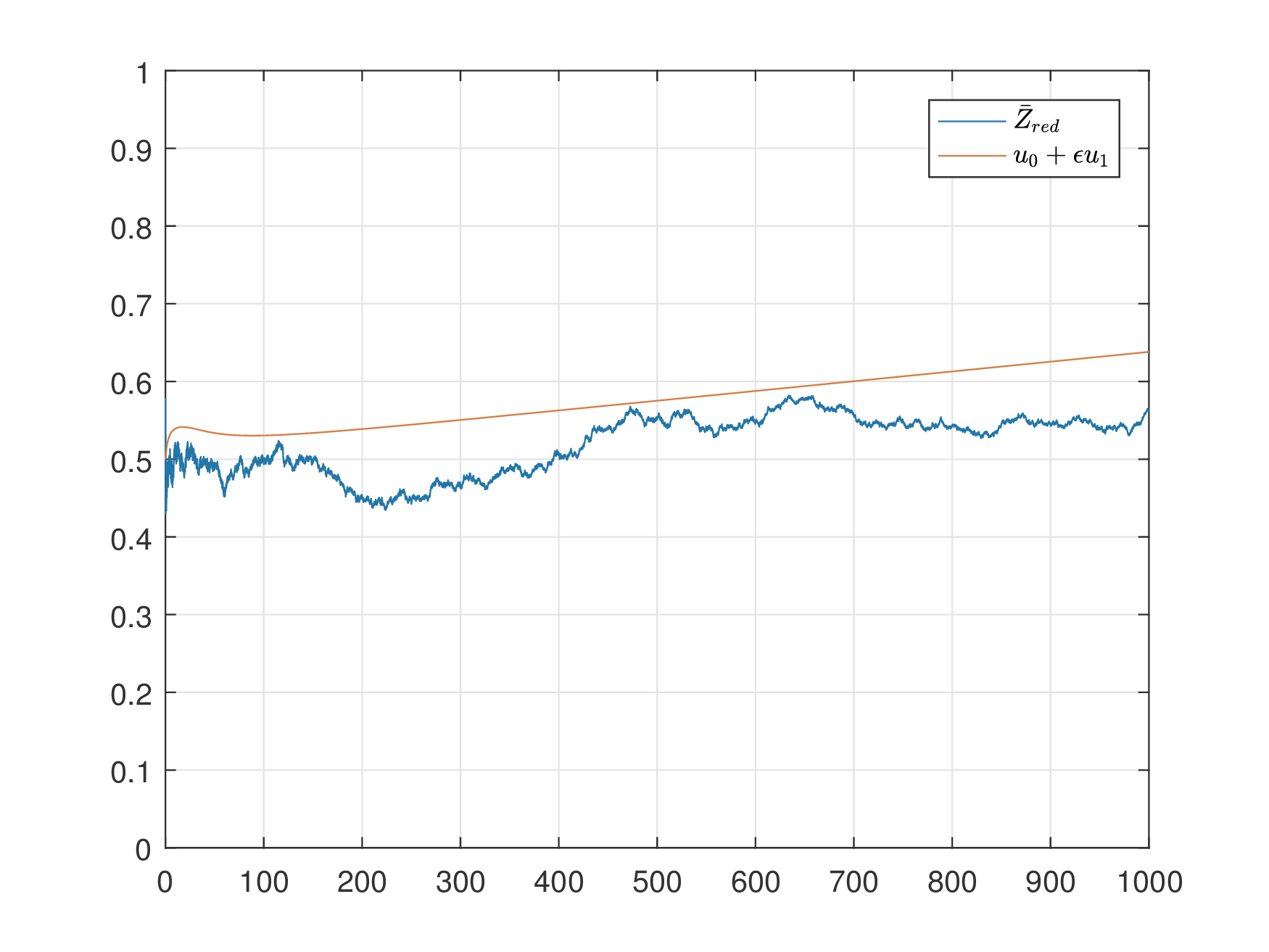}
            \caption{Approximation evolves to an $y$-bar state for $\hat{\epsilon}=0.1$.}\label{fig:approx minus}
        \end{subfigure}
        \hfill
        \begin{subfigure}[b]{0.475\textwidth}
            \centering
            \includegraphics[width=\textwidth]{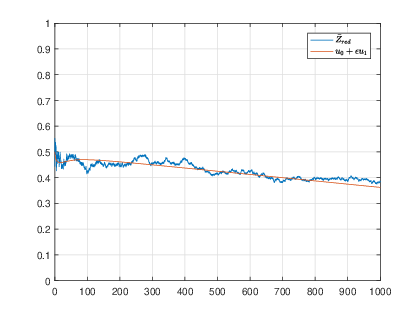}
            \caption{Approximation evolves to an $x$-bar state for $\hat{\epsilon}=-0.1$.}\label{fig:approx plus}
        \end{subfigure}
                \caption[ ]{Approximation follows the evolution of the $Z_{red}(t)$ for small $\epsilon$.}
        \label{fig:approx}
\end{figure}

%
Since the Monte Carlo simulation of \eqref{E:deltanot1} used zero initial conditions, we plot the approximation $\hat{u}^{\epsilon}(p,q,t)$ for $(p,q)=(0.1,0.1)$, close to the origin. With this choice, figures \ref{fig:approx minus} and \ref{fig:approx plus} show that in the perturbed system, the $\mathcal{O}(\epsilon)$ approximation to $u^\epsilon(\tau)=\mathbb{E}[Z_{red}(\tau)]$ evolves toward 0 or 1 depending on the sign of $\hat{\epsilon}:=\epsilon_0\epsilon.$

\begin{figure}[H]
	\centering
	\includegraphics[width=10cm, height=6cm]{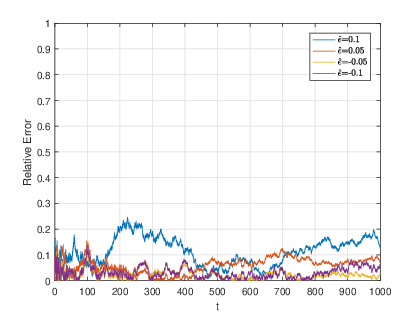}
	\caption{Relative Error $\frac{|\hat{u}^\epsilon(t)-\bar{Z}_{red}(t)|}{\bar{Z}_{red}(t)}$.}\label{fig:re}
\end{figure}

Using these simulations with $\nu$ scaled, we now explore how the intervals on which $\hat{u}^\epsilon$ serves as a good approximation to $\bar{Z}_{red}(t)$ depend on the perturbation parameter. Figure \ref{fig:re} shows the relative error (RE) given by
\[
RE = \frac{|\hat{u}^\epsilon(t)-\bar{Z}_{red}(t)|}{\bar{Z}_{red}(t)}.
\]

On some initial interval of time, the PDE approximation, $\hat{u}^\epsilon$ indeed serves as a close approximation to $\bar{Z}_{red}(t)$. Furthermore, as expected, the smaller the perturbation parameter $\epsilon$, the longer the approximation is valid.

\section{Concluding remarks and future directions}\label{S:Conclusion}

In this paper, we developed a finite dimensional SDE model that can be used to elucidate the dynamics of the 2D Navier-Stokes vorticity equation with noise. Monte Carlo simulation of the reduced model showed that the major qualatative property of the system, i.e. the dominant quasi-stationary state, can be determined from the model. In particular, as has been observed numerically and rigorously, the existence and attracting nature of these quasi-stationary states play an important role in the evolution of the stochastic Navier-Stokes vorticity equation. Specifically, the aspect ratio of the periodic domain, $D_\delta=[0,2\pi\delta]\times[0,2\pi]$, determines whether generic solutions evolve toward an $x$-bar state ($\delta>1$), a $y$-bar state ($\delta<1)$, or a dipole state ($\delta=1$).

Perturbation analysis then shows that the proposed reduced model can be viewed as a slow-fast system, Subsequent averaging and homogenization methods show the leading order behavior as the perturbation parameter $\delta\approx1$ goes to $\delta=1,$ in relation to how the viscocity parameter $\nu$ vanishes.

The numerical studies in \S\ref{S:sims} show that, on average, the system prefers to trend toward the appropriate quasi-stationary state as determined by $\delta$, see Figure \ref{fig: transition average Z}. However, one can see from the sample path plotted in Figure \ref{fig:sample_path}, individual sample paths do exhibit transitions between $x$-bar and $y$-bar states, as it has also been observed in \cite{BouchetSimonnet09}.

In regards to future directions, there are a number of interesting questions that one can ask and hope to answer. To begin with, the perturbation analysis of \S \ref{S:PDE} is formal and one would like to prove both well-posedness of (\ref{E:u0u1new}) and validity of the perturbation expansion.


In addition, the numerical studies of \S \ref{S:sims} suggest that while there are transitions at the individual sample path level, the system tends to converge to the preferred state depending on whether $\delta<1$ or $\delta>1$. One would like to make this mathematically rigorous. Furthermore, one could potentially use the reduced model of \S \ref{S:FourierSpaceRep} to build a related large deviations theory describing probabilities of the system being in, and exit times for leaving, one of the quasi-stationary states.

Furthermore, the form of the noise used in (\ref{E:noise}) and \eqref{E:xy}, with scaling factor $\sqrt{\nu}$, is common in the literature, see for example \cite{Glatt-Holtz} and the references therein. One of the important questions in the literature is the investigation of the convergence of the corresponding invariant measure as $\nu\to 0$ to that of  the 2D Euler equation, see \cite{kuksin_shirikyan_2012}. The support of the limiting invariant  measure (i.e. as $\nu\to 0$), is in general still an open question, see \cite{Glatt-Holtz}, and it has been resolved in special cases in  \cite{Bedrossian2} and \cite{Bedrossian1}. In our work, we study the long time behavior of the vorticity equation as $\nu\to0$ and in particular the selection mechanism for the dominant quasi-stationary states. It is reasonable to expect that there is a relation between the selection mechanism and the support of the limiting invariant measure. This is an intriguing question that is left for future work and it is beyond the scope of this work.

\section{Appendix}\label{Appendix1}
The complete list of operators in the Kolmogorov equation \eqref{E:kol} is given by
\begin{eqnarray}
\mathcal{L}_0u &=& -\frac{1}{5\nu_0}(p^2+q^2)\left(r\frac{\partial u}{\partial r} + s\frac{\partial u}{\partial s}\right)+2\nu_0\left(\sigma_5^2\frac{\partial^2 u}{\partial r^2}+\sigma_7^2  \frac{\partial^2 u}{\partial s^2}\right)\nonumber \\
\mathcal{L}_1u &=& -\frac{\epsilon_0}{100\nu_0}(51p^2-31q^2)\left(r\frac{\partial u}{\partial r} + s\frac{\partial u}{\partial s}\right) \nonumber \\
\mathcal{L}_2 u &=& -\frac{373}{1000\nu_0}(p^2+q^2)\left(r\frac{\partial u}{\partial r} + s\frac{\partial u}{\partial s}\right) \nonumber \\
\mathcal{L}_3 u &=&-\frac{\epsilon_0}{1000\nu_0}(379p^2-4109q^2)\left(r\frac{\partial u}{\partial r} + s\frac{\partial u}{\partial s}\right) \nonumber \\\nonumber \\
\mathcal{L}_{4}u &=& \frac{3}{40\nu_0}(r^2+s^2)\left(p\frac{\partial u}{\partial p}+q\frac{\partial u}{\partial q}\right)+2\nu_0\left(\sigma_1^2\frac{\partial^2}{\partial p^2}+\sigma_3^2\frac{\partial^2}{\partial q^2}\right) \nonumber\\
\mathcal{L}_{5}u &=&\frac{\epsilon_0}{\nu_0}(r^2+s^2)\left(\frac{27}{200}p\frac{\partial u}{\partial p}-\frac{3}{50}q\frac{\partial u}{\partial q}\right) \nonumber \\\nonumber \\
\mathcal{L}_6u &=& \frac{117}{4000\nu_0}(r^2+s^2)\left(p\frac{\partial u}{\partial p}+q\frac{\partial u}{\partial q}\right) \nonumber \\\nonumber \\
\mathcal{L}_7u &=&  -\frac{\epsilon_0}{\nu_0}(r^2+s^2)\left(\frac{123}{5000}p\frac{\partial u}{\partial p}+\frac{93}{20000}q\frac{\partial u}{\partial q}\right) \nonumber \\
\mathcal{L}_{8}u &=& -\frac{1}{2}(r-s)\left(q\frac{\partial u}{\partial p}-p\frac{\partial u}{\partial q}\right) \nonumber \\
\mathcal{L}_{9}u &=& \frac{\epsilon_0}{2}(r-s)\left(q\frac{\partial u}{\partial p}+p\frac{\partial u}{\partial q}\right) \label{E:multiscaleoperators} \\
\mathcal{L}_{10}u &=& \frac{1}{8}(r-s)\left(\frac{7}{2}q\frac{\partial u}{\partial p}-p\frac{\partial u}{\partial q}\right)\nonumber \\
\mathcal{L}_{11}u &=& \epsilon_0\frac{3}{8}q(r-s)\frac{\partial u}{\partial p}\nonumber 
\end{eqnarray}
\begin{eqnarray}
\mathcal{L}_{12}u &=& -\nu_0\left(p\frac{\partial u}{\partial p}+q\frac{\partial u}{\partial q}+2\left(r\frac{\partial u}{\partial r}+s\frac{\partial u}{\partial s}\right)\right)\nonumber \\
\mathcal{L}_{13}u &=&  \nu_0\epsilon_0\left(p\frac{\partial u}{\partial p}+r\frac{\partial u}{\partial r}+s\frac{\partial u}{\partial s}\right) \nonumber \\
\mathcal{L}_{14}u &=& -\nu_0\left(p\frac{\partial u}{\partial p}+r\frac{\partial u}{\partial r}+s\frac{\partial u}{\partial s}\right)\nonumber \\
\mathcal{L}_{15}u &=& \nu_0\epsilon_0\left(p\frac{\partial u}{\partial p}+r\frac{\partial u}{\partial r}+s\frac{\partial u}{\partial s}\right)\nonumber \\
\mathcal{L}_{16}u &=& -\epsilon_0pq\left(\frac{\partial u}{\partial r}-\frac{\partial u}{\partial s}\right) \nonumber \\
\mathcal{L}_{17}u &=& \frac{1}{2}pq\left(\frac{\partial u}{\partial r}-\frac{\partial u}{\partial s}\right) \nonumber \\
\mathcal{L}_{18}u &=& -\epsilon_0\frac{3}{8}pq\left(\frac{\partial u}{\partial r}-\frac{\partial u}{\partial s}\right)  \nonumber
\end{eqnarray}


\end{document}